\documentclass[11pt]{elsarticle}

\usepackage{amsmath,amsthm,bm}
\newtheorem{theorem}{Theorem}

\newtheorem{remark}{Remark}

\usepackage{caption}  
\usepackage{amsmath,amscd,wrapfig,amsfonts,amssymb,mathtools,mathrsfs}
\usepackage{booktabs,multirow,lscape,cool,bm}
\usepackage{url,parskip,caption}
\usepackage[shortlabels]{enumitem}
\usepackage[top=1.0in, bottom=1.0in, left=1.05in, right=1.05in]{geometry}

\DeclareMathOperator{\sign}{sign}

\DeclareMathOperator{\Ai}{Ai}
\DeclareMathOperator{\Bi}{Bi}
\DeclareMathOperator{\arccosh}{arccosh}

\makeatletter
\g@addto@macro\normalsize{%
  \setlength\abovedisplayskip{.7em}
  \setlength\belowdisplayskip{.7em}
  \setlength\abovedisplayshortskip{.7em}
  \setlength\belowdisplayshortskip{.7em}
}

\newcommand{\tcheb}[1]{t_{#1,k}^{\mbox{\tiny cheb}}}

\DeclareMathOperator{\res}{res}
\DeclareMathOperator{\diag}{diag}

\providecommand{\aiphase}{\alpha_{\mbox{\tiny ai}}}

\allowdisplaybreaks

\begin{document}

\begin{frontmatter}


\title
{
Airy Phase Functions
}

\begin{abstract}
It is well known that phase function methods allow for the numerical solution of a large class of oscillatory
second order linear ordinary differential equations in time independent of frequency.
Unfortunately,  these methods break down  in the commonly-occurring case in which the equation has 
turning points.  
Here, we resolve this difficulty by introducing a generalized phase function method 
designed for the case of   second order linear ordinary differential equations with turning points.
More explicitly, we prove the existence of a slowly-varying  ``Airy phase function'' 
that efficiently represents a basis in the space of solutions of such an equation,
and describe a numerical algorithm for calculating this Airy phase function.
The running time of our algorithm  is independent of the magnitude of the logarithmic 
derivatives of the equation's solutions, which is a measure of their rate of variation 
that generalizes the notion of frequency to functions which are rapidly varying
but not necessarily oscillatory.  Once the Airy phase function has been constructed, any reasonable initial or 
boundary value problem for the equation can be readily solved and, unlike step methods which output
the values of  a rapidly-varying solution on a sparse discretization grid that is insufficient for interpolation,   
the output of our scheme allows for the rapid evaluation of  the obtained solution at any  point in its domain.  
 We rigorously justify our approach by proving not only the existence of slowly-varying 
Airy phase functions, but also  the convergence of our numerical method.
Moreover, we present the results of extensive numerical experiments  demonstrating the efficacy 
of our algorithm.

\end{abstract}

\author{Richard Chow}
\address{\vskip -1em Department of Mathematics, University of Toronto}

\author{James Bremer}
\ead{bremer@math.toronto.edu}
\address{\vskip -1em Department of Mathematics, University of Toronto\vskip -2.5em}

\end{frontmatter}


\begin{section}{Introduction}

Second order linear ordinary differential equations frequently arise  in numerical and
scientific computations.  Many of the equations which are encountered are either of the form
\begin{equation}
y''(t) + \omega^2 q(t,\omega) y(t) = 0,\ \ \ \ \ a < t < b,
\label{introduction:ode}
\end{equation}
with $\omega$ a real-valued parameter  and $q$ a smooth function such that  
$q$ and its derivatives with respect to $t$ are bounded independent  of $\omega$,  or they can be easily 
put into  this form.      Several widely-used families of special functions,  
such as the Jacobi polynomials and the spheroidal wave functions, satisfy equations of this type.
They also arise in computations related to plasma physics \cite{Hazeltine,Davidson},  
Hamiltonian dynamics \cite{2018_Pritula} and cosmology \cite{Jerome,Agocs}, to name just a few 
representative applications.

For the sake of notational simplicity, we will generally suppress the dependence of $q$ on $\omega$.
This causes no harm because of our assumption that $q$ and its derivatives with respect to $t$ of all orders
are bounded independent of $\omega$.    Moreover, we will focus on  the commonly-occurring case in which
 $q$ has a single simple turning point in its domain $[a,b]$.  Without loss of generality, we will suppose that
 $a < 0 < b$,   and that $q(t) \sim t$ as $t \to 0$.  In particular, $q$ is strictly positive on $(0,b)$ and
strictly negative on $(a,0)$.
It is well known that, under these conditions, there exists a basis  $\{y_1,y_2\}$ in the space of solutions of (\ref{introduction:ode})
which can be asymptotically  approximated on various portions of the  interval $[a,b]$ as follows:
\begin{enumerate}[(a)]
\item\label{introduction:regimea} 
In any compact subinterval of  $(0,b]$, the estimates
\begin{equation}
y_1(t) = \frac{\cos\left(\alpha_0(t)\right)}{\sqrt{\alpha_0'(t)}} 
\left( 1 + \mathcal{O}\left(\frac{1}{\omega}\right)\right)
 \ \mbox{and}\ \,
y_2(t) = 
\frac{\sin\left(\alpha_0(t)\right)}{\sqrt{\alpha_0'(t)}} 
\left( 1 + \mathcal{O}\left(\frac{1}{\omega}\right)\right)
\label{introduction:uvtrig}
\end{equation}
hold uniformly with respect to $t$ as $\omega \to \infty$, where $\alpha_0$ is given by
\begin{equation}
\alpha_0(t) = \frac{\pi}{4} + \omega \int_0^t \sqrt{q(s)}\, ds.
\label{introduction:alpha0}
\end{equation}
\item\label{introduction:regimeb}
In any compact interval contained in  $[a,0)$, the estimates
\begin{equation}
y_1(t) = \frac{\exp\left(-\beta_0(t)\right)}{\sqrt{\beta_0'(t)}} 
\left( 1 + \mathcal{O}\left(\frac{1}{\omega}\right)\right)
 \ \mbox{and}\ \,
y_2(t) = \frac{\exp\left(\beta_0(t)\right)}{2 \sqrt{\beta_0'(t)}} 
\left( 1 + \mathcal{O}\left(\frac{1}{\omega}\right)\right)
\label{introduction:uvexp}
\end{equation}
hold uniformly with respect to $t$ as $\omega \to \infty$, where $\beta_0$ is defined via
\begin{equation}
\beta_0(t) =  \omega \int_0^t \sqrt{-q(s)}\, ds.
\label{introduction:beta0}
\end{equation}
%
\item\label{introduction:regimec}
In any compact subinterval of $[a,b]$, we have
\begin{equation}
y_1(t) = 
\frac{\Bi\left(\gamma_0(t)\right)} {\sqrt{\gamma_0'(t)}}
\left(1+ \mathcal{O}\left(\frac{1}{\omega}\right)\right)
\  \mbox{and} \ \,
y_2(t) = 
\frac{\Ai\left(\gamma_0(t)\right)} {\sqrt{\gamma_0'(t)}}
\left(1+ \mathcal{O}\left(\frac{1}{\omega}\right)\right)
\label{introduction:uvairy}
\end{equation}
uniformly with respect to $t$ as $\omega \to\infty$,  where $\gamma_0$ is defined by 
\begin{equation}
\gamma_0(t) = \begin{cases}
 \left(\frac{3}{2} \omega \int_{0}^t \sqrt{q(s)}\, ds \right)^{\frac{2}{3}} & \mbox{ if } t \geq 0\\
-\left(-\frac{3}{2}\omega \int_{0}^{t} \sqrt{-q(s)}\, ds \right)^{\frac{2}{3}} & \mbox{ if } t < 0
\end{cases}
\label{introduction:gamma0}
\end{equation}
and $\Ai$ and $\Bi$ refer to the solutions of Airy's differential equation
\begin{equation}
y''(t) + t y(t) = 0
\label{introduction:airyeq}
\end{equation}
such that
\begin{equation}
\begin{aligned}
\Ai(0) &= \frac{\sqrt{\pi}}{3^{\frac{2}{3}} \Gamma\left(\frac{2}{3}\right)}, 
 \ \ \Ai'(0) = \frac{\sqrt{\pi}}{3^{\frac{1}{3}}\Gamma\left(\frac{1}{3}\right)},\ \ \ 
\Bi(0) &= \frac{\sqrt{\pi}}{3^{\frac{1}{6}}\Gamma\left(\frac{2}{3}\right)} \ \ \ \mbox{and}\ \ \
\Bi'(0)= \frac{-3^{\frac{1}{6}}\sqrt{\pi}}{\Gamma\left(\frac{1}{3}\right)}.
\end{aligned}
\end{equation}
\end{enumerate}
%
We note that our definitions of $\Ai$ and $\Bi$ are
nonstandard.  It is more common for authors to work with the functions 
\begin{equation}
\frac{1}{\sqrt{\pi}}\, \Ai(-t) 
\ \ \ \mbox{and}\ \ \ 
\frac{1}{\sqrt{\pi}}\, \Bi(-t),
\end{equation}
which are solutions of 
\begin{equation}
y''(t) - t y(t) = 0.
\label{introduction:airybad}
\end{equation}
These are the definitions of the Airy functions which are used in \cite{DLMF} and \cite{Olver}, for example.
However, our convention is more suitable when treating
equations of the form (\ref{introduction:ode}) for which  $q(t) \sim t$ as $t \to 0$,
and  our choice of normalization ensures that the Wronskian of $\{\Ai(t),\, \Bi(t)\}$, and hence
also of the pair of asymptotic approximations appearing in (\ref{introduction:uvairy}),  is $1$.

It follows from (\ref{introduction:uvtrig}), (\ref{introduction:uvexp}) and
(\ref{introduction:uvairy}) that the solutions of an equation of the form (\ref{introduction:ode}) 
are rapidly varying  when $\omega$ is large.  Indeed,   
if $f_0(t) = y_1(t) + i y_2(t)$ and $j$ is a positive integer, then
\begin{equation}
\sup_{a\leq t \leq b} \left| \frac{ f_0^{(j)}(t) } {f_0(t)} \right|
= \mathcal{O}\left(\omega^j\right)
\ \mbox{as}\ \omega\to \infty.
\end{equation}
%
One obvious consequence  is that the cost of representing the solutions of (\ref{introduction:ode}) 
over the interval $[a,b]$ via polynomial or piecewise polynomial expansions grows linearly with $\omega$.
Since standard solvers for ordinary differential equations use such expansions to represent solutions either 
explicitly or implicitly, their running times also increase linearly with $\omega$.

In contrast to the rapidly-varying solutions of (\ref{introduction:ode}), $\alpha_0$, $\beta_0$ and $\gamma_0$ 
are slowly varying in the sense that whenever $f$ is equal to one of these functions, $[c,d]$ is a compact 
interval in the domain of definition of $f$ and $j$ is a positive integer, the quantity
\begin{equation}
\sup_{c\leq t \leq d} \left| \frac{ f^{(j)}(t) } {f(t)} \right|
\label{introduction:f}
\end{equation}
is bounded independent of $\omega$.  It follows that $\alpha_0$, $\beta_0$ and $\gamma_0$ can
be represented via polynomial expansions at a cost which is independent of $\omega$.
In particular, initial and boundary value problems for (\ref{introduction:ode}) can be solved with $\mathcal{O}\left(\omega^{-1}\right)$
relative accuracy in time independent of $\omega$ by forming polynomial expansions
of $\alpha_0$, $\beta_0$ and $\gamma_0$ and making use of the approximations 
(\ref{introduction:uvtrig}), (\ref{introduction:uvexp}) and (\ref{introduction:uvairy}).


It is well known that there exist higher order generalizations $\alpha_M$ and $\beta_M$ of the approximates 
$\alpha_0$ and $\beta_0$.  Indeed,  for each nonnegative integer $M$   and  compact subinterval $[c,d]$ of $(0,b]$, 
there exist smooth functions $\alpha:[c,d]\to\mathbb{R}$ and 
$\alpha_M:[c,d]\to\mathbb{R}$ such that
$\alpha_M$ is slowly varying in the aforementioned sense,
\begin{equation}
\left\{\frac{\cos\left(\alpha(t)\right)}{\sqrt{\alpha'(t)}},\ 
\frac{\sin\left(\alpha(t)\right)}{\sqrt{\alpha'(t)}} \right\}
\label{introduction:trigbasis}
\end{equation}
is a basis in the space of solutions of (\ref{introduction:ode}) given on the interval $[c,d]$ and
\begin{equation}
\alpha(t) = \alpha_M(t) \left(1 + \mathcal{O}\left(\frac{1}{\omega^{2(M+1)}}\right)\right)\ \mbox{as}\ \omega\to\infty.
\label{introduction:alphaestimate}
\end{equation}
Similarly,  for each compact subinterval $[c,d]$  of $[a,0)$ 
and  nonnegative integer $M$, there exist smooth functions
$\beta:[c,d]\to\mathbb{R}$ and $\beta_M:[c,d]\to\mathbb{R}$ such that 
$\beta_M$ is slowly varying,
\begin{equation}
\left\{\frac{\exp\left(-\beta(t)\right)}{\sqrt{\beta'(t)}},\ 
\frac{\exp\left(\beta(t)\right)}{2\sqrt{\beta'(t)}} \right\}
\label{introduction:expbasis}
\end{equation}
is a basis in the space of solutions of (\ref{introduction:ode}) given on the interval $[c,d]$ and
\begin{equation}
\beta(t) = \beta_M(t) \left(1 + \mathcal{O}\left(\frac{1}{\omega^{2(M+1)}}\right)\right)\ \mbox{as}\ \omega\to\infty.
\label{introduction:betaestimate}
\end{equation}
These estimates are typically derived by analyzing the Riccati equation
\begin{equation}
r'(t) + (r(t))^2 + \omega^2 q(t) = 0
\label{introduction:riccati}
\end{equation}
satisfied by the logarithmic derivatives of the solutions of (\ref{introduction:ode}); see, for instance,
\cite{StojimirovicBremer}.  We will refer to any $\alpha$ such that (\ref{introduction:trigbasis})
is a basis in the space of solutions of (\ref{introduction:ode}) on some subinterval of $[a,b]$ as a trigonometric phase function
for (\ref{introduction:ode}), and we call any $\beta$  such that 
(\ref{introduction:expbasis}) is a basis in the space of solutions of (\ref{introduction:ode}) on some subinterval
of $[a,b]$ an exponential phase function for (\ref{introduction:ode}).

While $\alpha_0$ and $\beta_0$ admit simple expressions 
that are easy to evaluate numerically, the same is not true of  $\alpha_M$ and $\beta_M$ when $M$ is large.
These functions are given by complicated expressions that involve repeated integration and differentiation,
with the consequence that neither symbolic nor numerical calculations allow for their effective evaluation.
Several recently introduced methods which apply to equations of the form (\ref{introduction:ode})
in the oscillatory regime (i.e., on regions where  $q$ is strictly positive)
overcome this difficulty by  numerically solving various nonlinear ordinary differential equations
to compute a trigonometric phase function $\alpha$ directly
rather than attempting to construct one of the approximates $\alpha_M$.
The first such scheme, that of \cite{BremerPhase},  considers the nonlinear differential equation
\begin{equation}
\omega^2 q(t)  -  (\alpha'(t))^2  + \frac{3}{4} \left(\frac{\alpha''(t)}{\alpha'(t)}\right)^2
-\frac{1}{2} \frac{\alpha'''(t)}{\alpha'(t)} = 0,
\label{introduction:kummer}
\end{equation}
which we refer to as Kummer's equation in  light of the article \cite{Kummer}.
The solutions of (\ref{introduction:kummer}) comprise the trigonometric phase functions for (\ref{introduction:ode}),
and  the method of \cite{BremerPhase} exploits the fact that  (\ref{introduction:alphaestimate})  implies  the existence of a 
solution of (\ref{introduction:kummer}) that  can be regarded as slowly varying for the purposes of  numerical computation.  
Indeed, by choosing $M$ to be large, we see that there exists a trigonometric phase function 
 which agrees to machine precision accuracy with a slowly-varying function even for modest values of $\omega$.
Of course, since almost all solutions of Kummer's equation are rapidly varying, some mechanism must be used to select 
this slowly-varying solution.  The algorithm of \cite{BremerPhase} employs a ''windowing scheme''
which involves smoothly deforming the coefficient $q$ into a constant over a portion of the solution domain.  

The paper \cite{StojimirovicBremer} introduces a faster and simpler algorithm
which appears to be the current state of the art  for solving equations of the form (\ref{introduction:ode})
in the oscillatory case.   It operates by numerically constructing a
slowly-varying solution of the Riccati equation (\ref{introduction:riccati}).
Once this has been done, a slowly-varying trigonometric phase function
$\alpha$ is constructed via the formula
\begin{equation}
r(t) = i \alpha'(t) - \frac{1}{2} \frac{\alpha''(t)}{\alpha'(t)}
\label{introduction:alphar}
\end{equation}
that connects the solutions of (\ref{introduction:riccati}) to  those of (\ref{introduction:kummer}).
To identify the desired  solution,
(\ref{introduction:riccati}) is discretized over a small interval $[a_0,b_0]$ in the solution domain via a Chebyshev spectral method
and Newton's method is applied to the resulting system of nonlinear algebraic equations.
The discretization grid is chosen to be dense enough to represent the desired slowly-varying
solution, but sparse enough that it does not discretize the rapidly-varying  ones.    
A proof that  the Newton iterates converge to a unique vector 
representing the desired solution of the Riccati equation provided that $\omega$ is sufficiently large is given in \cite{StojimirovicBremer}.
Once the values of this solution are known on the small interval $[a_0,b_0]$, an initial
value problem for (\ref{introduction:riccati}) is  solved in order to extend it 
to the entire solution domain and the desired slowly-varying
trigonometric phase function is then  constructed from $r$.

Both of the phase function methods of \cite{BremerPhase} and \cite{StojimirovicBremer}   run in time independent of $\omega$,
allow for the solution of essentially any reasonable initial or boundary value problem for (\ref{introduction:ode})
and  achieve accuracy on the order of the condition number of the problem being solved.
Moreover, unlike step methods which output the values of a rapidly-varying solution of (\ref{introduction:ode})
on a sparse grid of discretization nodes that does not suffice for interpolation, 
the output of \cite{BremerPhase} and \cite{StojimirovicBremer} can be 
used to  evaluate the obtained solution anywhere on the interval $[a,b]$ at a cost independent of $\omega$.

While it is shown in \cite{BremerPhase2} that 
trigonometric phase functions can represent the solutions of an equation of the form (\ref{introduction:ode}) near a turning point, 
they exhibit complicated behavior there, with the consequence that  the cost to calculate trigonometric phase functions
for second order linear ordinary differential equations with turning points is high.
Here, we address this difficulty by introducing a generalized phase function method
adapted to equations of the form (\ref{introduction:ode}).

We first show that, under our assumptions on $q$ and for any nonnegative integer $M$,  there exist a smooth function 
$\gamma:[a,b]\to\mathbb{R}$ and  a smooth slowly-varying function $\gamma_M:[a,b]\to\mathbb{R}$ such that 
\begin{equation}
\left\{\frac{\Bi\left(\gamma(t)\right)}{\sqrt{\gamma'(t)}},\ 
\frac{\Ai\left(\gamma(t)\right)}{\sqrt{\gamma'(t)}}\right\}
\label{introduction:airyphase}
\end{equation}
is a basis in the space of solutions of (\ref{introduction:ode}) and 
\begin{equation}
\gamma(t)= \gamma_M(t) 
\left(
1+\mathcal{O}\left(\frac{1}{\omega^{2(M+1)}}\right)
\right)
 \ \mbox{as}\ \omega\to\infty.
\label{introduction:gammaestimate}
\end{equation}
We will call any function $\gamma$ such that 
(\ref{introduction:airyphase}) is a basis in the space
of solutions of (\ref{introduction:ode}) an Airy phase function for (\ref{introduction:ode}).
While the estimates (\ref{introduction:alphaestimate}) and  (\ref{introduction:betaestimate})
follow easily from well-known results which can be found in many sources, the authors 
were unable to find a proof of (\ref{introduction:gammaestimate}) in the literature.
The standard asymptotic approximations of the solutions of (\ref{introduction:ode}),
which are discussed in  Chapter~11  of \cite{Olver}, among many other sources,
are of the form
\begin{equation}
\frac{z\left(\omega^{\frac{2}{3}}\xi(t)\right)}{\sqrt{\xi'(t)}}\, \sum_{j=0}^M \frac{A_j(\xi(t))}{\omega^{2j}}
+ 
\frac{z'\left(\omega^{\frac{2}{3}}\xi(t)\right)}{\omega^{\frac{4}{3}}\sqrt{\xi'(t)}}\, \sum_{j=0}^M \frac{B_j(\xi(t))}{\omega^{2j}},
\label{introduction:wkb}
\end{equation}
where  $\xi(t) = \omega^{-\frac{2}{3}} \gamma_0(t)$ 
%
and $z$ is taken to be  $\Ai$ or $\Bi$.

Airy phase functions satisfy the nonlinear differential equation
\begin{equation}
\omega^2 q(t)  -  \gamma(t) (\gamma'(t))^2  + \frac{3}{4} \left(\frac{\gamma''(t)}{\gamma'(t)}\right)^2
-\frac{1}{2} \frac{\gamma'''(t)}{\gamma'(t)} = 0,
\label{introduction:airykummer}
\end{equation}
which we will refer to as the Airy-Kummer equation,  and we go on to prove that they can be computed numerically by 
applying a method analogous to the algorithm of  \cite{StojimirovicBremer} to (\ref{introduction:airykummer}).
That is, we show that  if (\ref{introduction:airykummer}) is discretized over a small
interval $[-a_0,a_0]$ via a Chebyshev spectral method and Newton's method
is applied to the resulting system of nonlinear algebraic equations with the first-order approximate $\gamma_0$
used to form an initial guess, the iterates converge
to a vector representing the desired  slowly-varying Airy phase function $\gamma$, provided 
$a_0$ is sufficiently small and  $\omega$ is sufficiently large.
We use $\gamma_0$ as an initial guess because, unlike the higher order approximates $\gamma_M$,
it can be readily calculated.  
The analysis presented here  to establish the estimate (\ref{introduction:gammaestimate}) and prove
the convergence of the Newton iterations for (\ref{introduction:airykummer})
is significantly more involved that that of \cite{StojimirovicBremer} owing to the much more complicated structure 
of the Airy-Kummer equation (\ref{introduction:airykummer}) vis-\`a-vis the Riccati equation 
(\ref{introduction:riccati}).  

As  a last step in our algorithm, the function $\gamma$ is extended to the entire solution
domain $[a,b]$ of (\ref{introduction:ode}) by solving an initial value problem
for (\ref{introduction:airykummer}).    Here, we use
an adaptive Chebyshev spectral method, and the result is a collection of   piecewise
Chebyshev expansions representing the Airy phase function $\gamma$ and its first few derivatives.
Once these expansions have been obtained, any reasonable initial or boundary value problem for (\ref{introduction:ode}) 
can be easily solved.  The running time of our algorithm is independent of $\omega$
and the obtained solution of 
(\ref{introduction:ode}) can be readily evaluated at any point in the solution domain
at a cost which is independent of $\omega$.

The scheme of this paper can be easily combined with the phase function method of
\cite{StojimirovicBremer} by computing trigonometric or exponential phase functions in some regions
and Airy phase functions in others.  This is sometimes more efficient than using these algorithms 
separately.  Moreover, our method generalizes in a straightforward way to equations with multiple simple 
turning points: one simply constructs multiple Airy phase functions, one for each turning point, using the 
algorithm of this paper.   The situation regarding higher order turning points, however, is somewhat more complicated.
While numerical experiments show that a straightforward generalization of our numerical algorithm  works well 
in the setting of equations of the form (\ref{introduction:ode})
in which $q(t) \sim t^\sigma$ as $t\to 0$ with $\sigma > -2$, the analysis presented here
only applies when  $\sigma$ is $-1$, $0$ or $1$.    In this paper, we will focus on equations with simple turning points 
and we leave a treatment of the more general case for future work.

The remainder of this paper is structured as follows.    
 Section~\ref{section:preliminaries}
discusses  the necessary  mathematical and numerical preliminaries, and introduces the notation
used in the rest of the paper.  Our proof of the existence of slowly-varying Airy phase functions is given in Section~\ref{section:analysis}, while 
Section~\ref{section:newton} gives a proof that Newton's method converges when it is applied
to a spectral discretization of (\ref{introduction:airykummer}) on a small interval containing the turning point.
Section~\ref{section:algorithm} details our algorithm for the numerical computation of Airy phase functions,
and we present the results of numerical experiments conducted to assess its effectiveness in 
Section~\ref{section:experiments}.  We close with  a few brief remarks regarding this work and future directions 
for research in  Section~\ref{section:conclusion}.

\end{section}

\begin{section}{Preliminaries}
\label{section:preliminaries}

\begin{subsection}{Notation and conventions}
\label{section:preliminaries:notation}

We use $\partial_j K(t,s)$ for the partial derivative of a multivariate function
$K:\mathbb{R}^n\to\mathbb{R}$ with respect to its $j^{th}$ argument.
We denote the  Fr\'echet derivative of a map $F:X\to Y$ between Banach spaces  at the point $x$
by  $D_x F$.   The space of bounded linear functions $X \to Y$ is  $L(X,Y)$,
and we say  a function $F:\Omega \subset X \to Y$ given on an open subset $\Omega$ in $X$
is continuously differentiable  provided the map $\Omega \to L(X,Y)$ given by $x \mapsto D_xF$ 
is continuous.

We use  $\|f\|_\infty$ for the uniform  norm of the function $f$ over its domain of definition.
We let $C^1\left(\left[c,d\right]\right)$ be the Banach space of functions $f:[c,d]\to\mathbb{R}$ which have 
uniformly continuous derivatives.  When working with  $C^1\left(\left[c,d\right]\right)$,
we will use the nonstandard norm
\begin{equation}
\|f\|_\omega =\omega^2 \left\|f\right\|_\infty + \omega\left\|f'\right\|_\infty,
\label{preliminaries:wnorm}
\end{equation}
which is defined for any positive value of $\omega$.
We will use bold symbols for elements of the  Euclidean space $\mathbb{R}^k$,
and denote the $l^\infty\left(\mathbb{R}^k\right)$ norm of a vector
$\boldsymbol{v}$ by $\left\|\boldsymbol{v}\right\|_\infty$.  
We use $\diag\left(\boldsymbol{x}\right)$ for the $k\times  k$ diagonal matrix 
whose diagonal entries are the elements of the vector $\boldsymbol{x} \in \mathbb{R}^k$.
The Hadamard or pointwise product of the vectors $\boldsymbol{x}$ and 
$\boldsymbol{y}$ is denoted $\boldsymbol{x} \circ \boldsymbol{y}$.
Moreover,  if $\boldsymbol{v} = \left(\begin{array}{cccc} v_1 & v_2 & \cdots & v_k \end{array}\right)^\top$, 
then, for each integer $j$, we use $\boldsymbol{v}^{\circ j}$ for the vector
\begin{equation}
\left(\begin{array}{cccc} v_1^j & v_2^j & \cdots & v_k^j \end{array}\right)^\top.
\end{equation}

We denote the Chebyshev polynomial of degree $j$ by $T_j$, 
and use $\tcheb{1}, \tcheb{2}, \ldots, \tcheb{k}$  for the nodes of the  $k$-point Chebyshev extremal grid on $[-1,1]$, 
which are given by  the formula
\begin{equation}
\tcheb{i} =  \cos\left(\pi \frac{k-i}{k-1}\right).
\end{equation}
The $k\times k$ identity matrix is $\mathscr{I}_k$, 
and we use $\mathscr{D}_k$ for the $k\times k$ spectral differentiation matrix
which takes the vector
\begin{equation}
\left(
\begin{array}{cccccc}
p\left(\tcheb{1}\right) & p\left(\tcheb{2}\right) & \cdots & p\left(\tcheb{k}\right)
\end{array}
\right)^\top
\label{preliminaries:pvals}
\end{equation}
of values of a polynomial $p$ of degree less than $k$  at the Chebyshev nodes to the vector
\begin{equation}
\left(
\begin{array}{cccccc}
p'\left(\tcheb{1}\right) & p'\left(\tcheb{2}\right) & \cdots & p'\left(\tcheb{k}\right)
\end{array}
\right)^\top
\end{equation}
of the values of its derivatives at the same nodes.  

\end{subsection}

\begin{subsection}{The Newton-Kantorovich theorem}

In \cite{Kantorovich},  Kantorovich generalized Newton's method to the case of maps between Banach spaces and
gave conditions for its convergence.  Here, we state a simplified version of the 
Newton-Kantorovich theorem that can be found as Theorem~7.7-4 in Section~7.7 of \cite{Ciarlet}.   
Note that we have corrected a minor typo in condition (5) of the theorem.
\begin{theorem}
\label{theorem:nk}
Suppose that $\Omega$ is an open subset of the Banach space $X$, 
that $Y$ is a Banach space and that $F: \Omega \subset X \to Y$ is continuously differentiable.
Suppose also that there exist a point $x_0 \in \Omega$ and
constants $\lambda$ and $\eta$  such that
\begin{enumerate}

\item
$\begin{aligned}
D_{x_0} F
\end{aligned}$
admits an inverse $\left(D_{x_0}F\right)^{-1} \in L(Y,X)$,

\item 
$\begin{aligned}
B_\eta\left(x_0\right) \subset \Omega,
\end{aligned}$

\item
$\begin{aligned}
0 < \lambda < \frac{\eta}{2},
\end{aligned}$

\item
$\begin{aligned}
\left\| \left(D_{x_0}F\right)^{-1} F\left(x_0\right) \right\| \leq \lambda \ \ \ \mbox{and}
\end{aligned}$

\item
$\begin{aligned}
\left\|\left(D_{x_0}F\right)^{-1} \left(D_xF - D_yF \right) \right\| \leq \frac{1}{\eta} \left\|x-y\right\|
\  \mbox{for all}  \  x,y \in B_\eta(x_0).
\end{aligned}$
\end{enumerate}
Then, $D_xF$ has a bounded inverse $D_xF^{-1} \in L(Y,X)$ for each $x \in B_\eta(x_0)$, 
the sequence $\{x_n\}$ defined by
\begin{equation}
x_{n+1} = x_n - \left(D_{x_n}F\right)^{-1} F(x_n)
\end{equation}
is contained in the open ball
\begin{equation}
 B_{\eta^-}(x_0), \ \ \ \mbox{where} \ \ \ \eta^- = \eta \left(1 - \sqrt{1-\frac{2\lambda}{\eta}}\right) \leq \eta
\end{equation}
and $\{x_n\}$ converges to a zero $x^*$ of $F$.
Moreover,  $x^*$ is the only  zero of $F$ in the ball $B_\eta(x_0)$
and, for each $n \geq 0$, we have
\begin{equation}
\left\| x_n - x^* \right\| \leq \frac{\eta}{2^n} \left( \frac{\eta^-}{\eta}\right)^{2^n}.
\label{preliminaries:nkconclusion}
\end{equation}

\end{theorem}


\end{subsection}

\begin{subsection}{A Green's function for an auxiliary differential equation related to the Airy functions}
\label{section:preliminaries:green}

In Section~\ref{section:analysis}, we use a Green's function $G(t,s)$ for
the differential equation
\begin{equation}
z'''(t) + 4 t z'(t) + 2 z(t) = 0
\label{preliminaries:airyappell}
\end{equation}
in order to form an integral operator to which we then apply a contraction mapping argument.
It is critical that $G(t,s)$ and its first order partial
derivatives  $\partial_1 G(t,s)$ and $\partial_2 G(t,s)$  be bounded for all
$(t,s) \in \mathbb{R}\times\mathbb{R}$.
In addition, it will be necessary to choose a Green's function  which  enforces the condition $z(0) =0$.

Fortunately, because the solutions of  (\ref{preliminaries:airyappell}) are explicitly known, 
it is straightforward to construct a Green's function with the desired properties.
Indeed,
\begin{equation}
\varphi_1(t) = 2^{-\frac{1}{3}} \Ai^2(t), \quad
\varphi_2(t) = 2^{-\frac{1}{3}} \Bi^2(t) \quad\mbox{and}\quad
\varphi_3(t) = 2^{-\frac{1}{3}} \Ai(t)\Bi(t) 
\end{equation}
is a triple of solutions of (\ref{preliminaries:airyappell})  whose Wronskian is $1$,
and  the method of  variations of parameters can be used to show that
\begin{equation}
G_0(t,s) =  \begin{cases}
-\varphi_2(t) \left(\varphi_1(s)\varphi_3'(s) - \varphi_1'(s)\varphi_3(s) \right) +  
\varphi_3(t) \left(\varphi_1(s)\varphi_2'(s)-\varphi_1'(s)\varphi_2(s)\right) 
& \mbox{if}\ s\leq t\\
-\varphi_1(t) \left(\varphi_2(s)\varphi_3'(s) - \varphi_2'(s)\varphi_3(s) \right) & \mbox{if}\ s\geq t
\end{cases}
\end{equation}
is a Green's function for (\ref{preliminaries:airyappell}).  
The well-known asymptotic estimates
\begin{equation}
\begin{aligned}
\Ai(t) \sim t^{-\frac{1}{4}} \cos\left(\frac{\pi}{4} - \frac{2}{3} t^\frac{3}{2}\right), \ \ \ \ \ \ \ \,
\Bi(t) \sim t^{-\frac{1}{4}} \sin\left(\frac{\pi}{4} - \frac{2}{3} t^\frac{3}{2}\right)
\ \ \ \mbox{as} \ \ t \to\infty
\end{aligned}
\label{preliminaries:airyasym1}
\end{equation}
and
\begin{equation}
\begin{aligned}
\Ai(t) \sim \frac{1}{2} t^{-\frac{1}{4}} \exp\left(- \frac{2}{3} (-t)^\frac{3}{2}\right),\ \ 
\ \ \Bi(t) \sim t^{-\frac{1}{4}} \exp\left( \frac{2}{3} (-t)^\frac{3}{2}\right)
\ \,  \ \ \mbox{as} \ \ t \to-\infty
\end{aligned}
\label{preliminaries:airyasym2}
\end{equation}
imply that  $G_0(t,s)$ is bounded on $\mathbb{R}\times\mathbb{R}$.  Moreover, the estimates
(\ref{preliminaries:airyasym1}) and (\ref{preliminaries:airyasym2}) can be differentiated
in order to see that the partial derivatives $\partial_1 G_0(t,s)$ 
and $\partial_2 G_0(t,s)$ are also bounded for all $(t,s) \in \mathbb{R}\times\mathbb{R}$.
Since $\Ai(0) \neq 0$, we can now define 
\begin{equation}
G(t,s) = G_0(t,s) - \frac{\Ai^2(t)}{\Ai^2(0)} G_0(0,s)
\label{preliminaries:G}
\end{equation}
so that 
\begin{equation}
z(t) = \int_c^d G(t,s) f(s)\, ds
\end{equation}
is a solution of the inhomogeneous equation
\begin{equation}
z'''(t) + 4 t z'(t) + 2 z(t)  = f(t), \quad c < t < d,
\end{equation}
such that $z(0) = 0$.    That $G(t,s)$ and its partial derivatives
$\partial_1 G(t,s)$ and $\partial_2 G(t,s)$ are bounded for all $(t,s) \in \mathbb{R}\times\mathbb{R}$
is clear given the properties of $G_0$ and $\Ai$.

\end{subsection}

\begin{subsection}{Generalized phase functions}
\label{section:liouville}

Trigonometric and exponential phase functions are particular cases of a more general construction
which we now discuss.  It can be easily seen that the transformation
\begin{equation}
y(t) = \frac{z(\lambda(t))}{\sqrt{\lambda'(t)}}
\label{liouville:forwardtransform}
\end{equation}
takes the solutions of  
\begin{equation}
z''(t) + \widetilde{q}(t) z(t) = 0
\label{liouville:eq2}
\end{equation}
to those of (\ref{introduction:ode}) provided $\lambda$ satisfies the nonlinear ordinary differential equation
\begin{equation}
\omega^2 q(t)  -  \widetilde{q}(\lambda(t)) (\lambda'(t))^2  + \frac{3}{4} \left(\frac{\lambda''(t)}{\lambda'(t)}\right)^2
-\frac{1}{2} \frac{\lambda'''(t)}{\lambda'(t)} = 0.
\label{liouville:genkummer}
\end{equation}
We will call (\ref{liouville:genkummer}) the generlized Kummer equation and refer to its solutions as 
generalized phase functions
for (\ref{introduction:ode}).
Their utility lies in the fact that the combination of 
the solutions of (\ref{liouville:eq2}) and the transformation $\lambda$ can be simpler in some fashion
 than  the solutions of  (\ref{introduction:ode}).    

When $q$ is strictly positive,  one usually chooses $\widetilde{q}$ to be $1$.
In this event, (\ref{liouville:genkummer}) becomes Kummer's equation (\ref{introduction:kummer})
and (\ref{liouville:eq2}) is
\begin{equation}
z''(t) + z(t) = 0.
\label{liouville:trig}
\end{equation}
Since $\left\{\cos(t),\sin(t)\right\}$ is a basis in the space of solutions of 
(\ref{liouville:trig}),  it follows that if $\alpha$ satisfies (\ref{introduction:kummer}),
then (\ref{introduction:trigbasis}) is a basis in the space of solutions of (\ref{introduction:ode}).  
That is, the solutions  of Kummer's equation are trigonometric phase functions.
When $q$ is strictly negative, one typically takes $\widetilde{q}=-1$ so that 
(\ref{liouville:genkummer}) becomes
\begin{equation}
\omega^2 q(t)  +  (\beta'(t))^2  + \frac{3}{4} \left(\frac{\beta''(t)}{\beta'(t)}\right)^2
-\frac{1}{2} \frac{\beta'''(t)}{\beta'(t)} = 0.
\label{liouville:kummerexp}
\end{equation}
Since, in this case, $\{\exp(-t),\exp(t)/2\}$ is a basis in the space of solutions of (\ref{liouville:eq2}),
(\ref{introduction:expbasis}) is a basis in the space of solutions of (\ref{introduction:ode})
whenever $\beta$ satisfies (\ref{liouville:kummerexp}).
That is, the solutions of (\ref{liouville:kummerexp}) are exponential phase functions for (\ref{introduction:ode}).

In this paper, our interest is in the case
in which the coefficient $q$ in (\ref{introduction:ode}) has a simple turning point.
In this event, it is more appropriate to take $\widetilde{q}(t)=t$ so that
(\ref{liouville:eq2}) becomes Airy's equation (\ref{introduction:airyeq})
and  (\ref{liouville:genkummer})  becomes the Airy-Kummer equation (\ref{introduction:airykummer}).
Since the Airy functions $\Ai(t)$ and $\Bi(t)$ constitute a basis in the space of solutions of (\ref{introduction:airyeq}),
the solutions of (\ref{introduction:airykummer})  are Airy phase functions for (\ref{introduction:ode}).

\end{subsection}

\end{section}

\begin{section}{Existence of slowly-varying Airy phase functions}
\label{section:analysis}

In this section, we prove the estimate (\ref{introduction:gammaestimate}).
Our argument is divided into three parts.   In Subsection~\ref{section:analysis:formal},
we use the method of matched asymptotic expansions 
to  show that there exists a sequence $u_0,u_1,\ldots,u_M$ of slowly-varying functions
such that when the function $\gamma_M$ defined via
\begin{equation}
\gamma_{M}(t) =  
\omega^{\frac{2}{3}}
\left(   u_0(t) + \frac{u_1(t)}{\omega^2} + \cdots + 
\frac{u_{M}(t)}{\omega^{2M}}\right)
\label{analysis:gammaM}
\end{equation}
is inserted into the Airy-Kummer equation, the residual
\begin{equation}
R_M(t) = \omega^2 q(t)  -  \gamma_M(t) (\gamma_M'(t))^2  + \frac{3}{4} \left(\frac{\gamma_M''(t)}{\gamma_M'(t)}\right)^2
-\frac{1}{2} \frac{\gamma_M'''(t)}{\gamma_M'(t)}
\label{analysis:residual}
\end{equation}
is on the order of $\omega^{-2M}$.   It will emerge that, for all sufficiently large $\omega$,
$\gamma_M'$ is positive on $[a,b]$ and $\gamma_M$ has a single zero $t_0$ in the interval $[a,b]$.
In Subsection~\ref{section:analysis:inteq}, we represent a solution $\gamma$
of the Airy-Kummer equation in the form
\begin{equation}
\gamma(t) = \int_{t_0}^t \gamma_M'(s)\exp\left(2\delta(s)\right)\, ds.
\label{analysis:gamma}
\end{equation}
and derive an integral equation for the function $\delta$.  
Finally, in Subsection~\ref{section:analysis:contraction},
we use a contraction mapping argument to show the existence of a solution of this integral equation
whose $L^\infty\left([a,b]\right)$ norm is on the order of  $\omega^{-2(M+1)}$.  
It follows by  applying the integral mean value theorem to  (\ref{analysis:gamma}) that, for each
$t \in [a,b]$, there exists a point $\xi$ in $[a,b]$ such that
\begin{equation}
\begin{aligned}
\gamma(t) = \left( \int_{t_0}^t \gamma_M'(s)\ ds \right) \exp(2\delta(\xi))
          = \gamma_M(t) \left(1 + 2 \delta(\xi) + (\delta(\xi)^2 + \cdots \right).
\end{aligned}
\label{analysis:imvt}
\end{equation}
%
The estimate (\ref{introduction:gammaestimate}) follows immediately
from (\ref{analysis:imvt}) and the fact that 
the  $L^\infty\left([a,b]\right)$ norm of $\delta$ is on the order of  $\omega^{-2(M+1)}$.  

\begin{subsection}{Formal asymptotic expansion}
\label{section:analysis:formal}

It is convenient to define the function
\begin{equation}
Q(t) = 
\int_0^t \sqrt{\left|q(s)\right|}\, ds.
\label{analysis:Q}
\end{equation}
By assumption,  $q$ is a smooth function such that $q(t)\sim t$ as $t\to 0$, 
and the only zero of $q$ in the interval $[a,b]$ occurs at the point $0$. 
These properties of $q$ together with  L'H\^opital's rule can be used to show  that  $Q$ and $Q'$ have the forms
\begin{equation}
Q(t) = \frac{2}{3} \sign(t) \left|t\right|^{\frac{3}{2}} Q_0(t)
\quad\mbox{and}\quad Q'(t) =  \sqrt{|t|} Q_1(t),
\label{analysis:Qform}
\end{equation}
where $Q_0$ and $Q_1$ are smooth and positive on $[a,b]$.
We now plug the ansatz (\ref{analysis:gammaM}) into the Airy-Kummer equation (\ref{introduction:airykummer}) 
and write the resulting expression as 
\begin{equation}
\begin{aligned}
\omega^2 q(t) - \omega^2 \left(\sum_{j=0}^M \frac{u_j(t)}{ \omega^{2j}} \right) \left(\sum_{j=0}^M \frac{u_j'(t)}{\omega^{2j}}\right)^2
+ \sum_{j=0}^\infty \frac{S_j(u_0(t),u_1(t),\ldots,u_{j}(t))}{\omega^{2j}} = 0,
\end{aligned}
\label{analysis:eq}
\end{equation}
where $S_j(u_0(t),\ldots,u_j(t))$ is the coefficient of $\omega^{-2j}$ in the expansion
of the Schwarzian derivative
\begin{equation}
\frac{3}{4} \left(\frac{\gamma_M''(t)}{\gamma_M'(t)}\right)^2
-\frac{1}{2} \frac{\gamma_M'''(t)}{\gamma_M'(t)}
\end{equation}
in inverse powers of $\omega$.  By grouping terms in  (\ref{analysis:eq}) we obtain 
an expression of the form
\begin{equation}
\begin{aligned}
&\omega^2 F_{-1}(Q(t),u_0(t))  + F_0(Q(t),u_0(t),u_1(t)) + \frac{1}{\omega^2} F_1(Q(t),u_0(t),u_1(t),u_2(t)) + \cdots +\\
&\quad\frac{1}{\omega^{2(M-1)}} F_{M-1}(Q(t),u_0(t),\ldots,u_{M}(t))  + \mathcal{O}\left(\frac{1}{\omega^{2M}}\right)
 = 0,
\end{aligned}
\label{analysis:expr}
\end{equation}
where
\begin{equation}
F_{-1}\left(Q(t),u_0(t)\right) = \sign(t)\left(Q'(t)\right)^2 - \left(u_0'(t)\right)^2\left(u_0(t)  \right)
\label{analysis:F2}
\end{equation}
and each subsequent $F_k$ takes the form 
\begin{equation}
F_k(Q(t),u_0(t),u_1(t),\ldots,u_{k}(t),u_{k+1}(t)) = 
-(u_0'(t))^2 u_{k+1}(t) -2 u_0(t) u_0'(t) u_{k+1}'(t) + \res_k(t)
\end{equation}
with
\begin{equation}
\res_k(t) = S_k(u_0(t),\ldots,u_{k}(t))\ - \sum\limits_{\substack{i_1+i_2+i_3=k\\0 \leq i_1,i_2,i_3 < k}} u_{i_1}(t) u_{i_2}'(t) u_{i_3}'(t).
\end{equation}
The nonlinear differential equation $F_{-1}(Q(t),u_0(t))=0$ admits the 
solution
\begin{equation} 
u_0(t) = 
\sign(t) \left(  \frac{3}{2} \sign(t) Q(t) \right)^{\frac{2}{3}}.
\label{analysis:u0}
\end{equation}
From (\ref{analysis:Qform}), we see that 
\begin{equation}
u_0(t) = t \left(Q_0(t)\right)^{\frac{2}{3}} \quad\mbox{and}\quad
u_0'(t) =  \frac{Q_1(t)}{\left(Q_0(t)\right)^{\frac{1}{3}}}
\end{equation}
so that $u_0$ is smooth, $u_0'$ is nonzero on $[a,b]$ and the only zero of $u_0$ on $[a,b]$
is the simple zero at $0$.
With our choice of $u_0$, $F_k$  becomes
\begin{equation}
F_k(Q(t),u_0(t),u_1(t),\ldots,u_{k}(t),u_{k+1}(t)) = 
-2^{\frac{2}{3}}
 Q'(t) \frac{d}{dt}\left[   
\left(3Q(t) \right)^{\frac{1}{3}}  
u_{k+1}(t)
\right] +\res_{k}(t).
\label{analysis:Fk}
\end{equation}
We now view $F_0(Q(t),u_0(t),u_1(t))=0$ as a  first order linear ordinary  differential equation defining $u_1$.  
The general solution is
\begin{equation}
u_1(t) = 
\left(3Q(t,\omega)\right)^{-\frac{1}{3}}
\left(C + 
 2^{-\frac{2}{3}} \,  \int_0^t  \frac{\res_0(s)}{Q'(s)}\, ds\right)
\label{analysis:gamma1}
\end{equation}
with $C$ an arbitrary constant.   We choose $C = 0$ so that $u_1$ is equal to
\begin{equation}
\frac{ 2^{-\frac{2}{3}}}{\sign(t) \sqrt{|t|} (3Q_0(t) )^{\frac{1}{3}}} \int_0^t \frac{\res_0(s) }{\sign(s) \sqrt{|s|}\sqrt{q_0(s)} }\, ds.
\end{equation}
It is easy to see that 
\begin{equation}
S_0(u_0(t)) = \frac{3}{4}\left(\frac{u_0''(t)}{u_0'(t)}\right)^2 - \frac{1}{2} \frac{u_0'''(t)}{u_0'(t)},
\end{equation}
which is smooth on $[a,b]$ since $u_0'(t)$ is nonzero there.  It follows that $\res_0(t)$ is also smooth.
In particular, $u_1(t)$ is of the form
\begin{equation}
\frac{f(t)}{k(t)} \int_0^t \frac{g(s)}{k(s)}\, ds
\label{analysis:smoothint}
\end{equation}
with $f$ and $g$ smooth functions and $k(s) = \sign(s) \sqrt{|s|}$. 
It is straightforward, but somewhat tedious, to use
L'H\^opital's rule to show that functions of the form (\ref{analysis:smoothint})  are smooth.  
We note that any other choice of $C$ would result in $u_1$ being singular at $0$.

We next observe that all of the $S_j$ are smooth when $\omega$ is sufficiently large.  
To see this, we first write 
\begin{equation}
\begin{aligned}
\frac{1}{\gamma_M'(t)} 
= 
\frac{1}{\omega^2 u_0'(t)} \left(1 + \frac{1}{\omega^2} \frac{u_1'(t)}{u_0'(t)} + \cdots + \frac{1}{\omega^{2M}}\frac{u_M'(t) }{u_0'(t)}\right)^{-1},
\end{aligned}
\end{equation}
which is allowable since $u_0'$ is positive on $[a,b]$.  Assuming $\omega$ is large enough,
the expression appearing inside the parentheses can be expanded in a  convergent Neumann series in order to obtain
formulas for the $S_j$ which show that they are smooth.
It follows that $\res_1(t)$ is smooth on $[a,b]$, and  we can  apply the preceding argument to show that there exists
a smooth solution $u_2$ of the first order linear ordinary differential equation
$F_1(Q(t),u_0(t),u_1(t),u_2(t))=0$ provided, of course, that $\omega$ is sufficiently large.  Continuing in this fashion gives us 
a sequence of smooth functions $u_0,u_1,\ldots,u_M$ such that the residual
(\ref{analysis:residual}) is on the order of $\omega^{-2M}$.
By construction, the functions $u_j$ are independent of $\omega$ when $q$ is independent of $\omega$
and slowly-varying in the sense discussed in the introduction when $q$ depends on $\omega$
but its derivatives are bounded independent of $\omega$.

We close this subsection by noting that, for all sufficiently large $\omega$, 
$\gamma_M'$ is nonzero on $[a,b]$,
 $\gamma_M$ has exactly one zero $t_0$ in $[a,b]$ and  $|t_0| = \mathcal{O}\left(\omega^{-2}\right)$ as $\omega \to \infty$.
This follows easily from the definition (\ref{analysis:gammaM}) of $\gamma_M$ and the observations about $u_0$ made above.


\begin{remark}
The method used here to construct a
formal asymptotic expansion representing a solution of  the Airy-Kummer equation (\ref{introduction:airykummer}) can be applied
to the  generalized Kummer equation (\ref{liouville:genkummer}) 
when $\widetilde{q}(t) = t^\sigma$, $q(t) = t^\sigma q_0(t)$ with $q_0$ a smooth positive function  and $\sigma$ equal 
to $-1$, $0$ or $1$.   Interestingly, though, it fails for all other values of $\sigma$.  Inserting 
\begin{equation}
\gamma(t) = \omega^{\frac{2}{2+\sigma}} \left(u_0(t) + \frac{u_1(t)}{\omega^2} + \cdots + \frac{u_{M}(t)}{\omega^{2M}}\right)
\end{equation}
into (\ref{liouville:genkummer}) yields a sequence of differential equations defining the $u_j$,
the first of which is the nonlinear equation
\begin{equation}
\sign(t) (Q'(t))^2 - \left(u_0'(t)\right)^2\left(u_0(t)  \right)^\sigma.
\label{analysis:remark:Qp}
\end{equation}
For all real-valued $\sigma > -2$, (\ref{analysis:remark:Qp}) admits the smooth solution
\begin{equation}
u_0(t) = 
\sign(t)  \left( \frac{2+\sigma}{2} \sign(t) \int_0^t \sqrt{\left|q(s)\right|}\, ds  \right)^{\frac{2}{2+\sigma}}.
\label{analysis:remark:u0}
\end{equation}
%
However, the subsequent equation which defines $u_1$ is
%
\begin{equation}
2^{\frac{2}{2+\sigma}} Q'(t) \frac{d}{dt} \left[
\left((\sigma+2) Q(t)\right)^{\frac{\sigma}{2+\sigma}} u_1(t) 
\right] = \res_0(t) :=
\frac{3}{4} \left(\frac{\gamma_0''(t)}{\gamma_0'(t)}\right)^2
-\frac{1}{2} \left(\frac{\gamma_0'''(t)}{\gamma_0'(t)}\right),
\end{equation}
and its general solution is 
\begin{equation}
u_1(t) = \frac{2^{-\frac{2}{2+\sigma}}}{\left((\sigma+2) Q(t)\right)^{\frac{\sigma}{2+\sigma}}} \left(C + \int_0^t \frac{\res_0(s) }{Q'(s)}\right),
\label{analysis:remark:u1}
\end{equation}
which has the  form
\begin{equation}
\frac{f(t)}{k(t)}
\left( C + 
\int_0^t \frac{g(s)}{k(s)}\, ds
\right)
\label{analysis:remark:u1form}
\end{equation}
with $k(t) \sim |t|^{\frac{\sigma}{2}}$ as $t\to 0$.
It can be readily seen that there is no choice of $C$ which makes  (\ref{analysis:remark:u1form})
smooth for all smooth $f$ and $g$ unless $\sigma$ is one of the special values  $-1$, $0$ or $1$.
\end{remark}

\end{subsection}

%
%
\begin{subsection}{Integral equation formulation}
\label{section:analysis:inteq}

We will now reformulate an integral equation for a function $\delta$ such that 
(\ref{analysis:gamma})
%
is a solution of the Airy-Kummer equation.  To that end, we first observe that 
if the function $w$ defined via
\begin{equation}
w(t) = \frac{1}{2} \log\left(\gamma'(t)\right)
\label{inteq:1}
\end{equation}
solves
\begin{equation}
w''(t) - (w'(t))^2 +  \exp(4 w(t)) \int_{t_0}^t \exp(2 w(s))\, ds- \omega^2 q(t) = 0,
\label{inteq:2}
\end{equation}
then the function (\ref{analysis:gamma}) satisfies (\ref{introduction:airykummer}).
We recall that  $t_0$ is chosen to be the sole zero of $\gamma_M$ in the interval
$[a,b]$.  Letting  
\begin{equation}
w(t) = w_M(t) + \delta(t)
\label{inteq:3}
\end{equation}
in (\ref{inteq:2}), where 
\begin{equation}
w_M(t) = \frac{1}{2} \log(\gamma_M'(t)),
\label{inteq:4}
\end{equation}
yields
\begin{equation}
\begin{aligned}
&\delta''(t) - 2w_M'(t)\delta'(t) - \left(\delta'(t) \right)^2 + \exp(4 w_M(t) + 4\delta(t))
\int_{t_0}^t \exp(2 w_M(s))\exp( 2\delta(s))\, ds  \\
&\quad =  -w_M''(t)+(w_M'(t))^2 + \omega^2 q(t).
\end{aligned}
\label{inteq:5}
\end{equation}
Subtracting $\exp(4 w_M(t)) \int_{t_0}^t \exp(2 w_M(s))\, ds$ from both sides of (\ref{inteq:5})
gives us 
\begin{equation}
\begin{aligned}
&\delta''(t) - 2w_M'(t)\delta'(t) - \left(\delta'(t) \right)^2 
+\exp(4 w_M(t))\exp( 4\delta(t)) \int_{t_0}^t \exp(2 w_M(s))\exp( 2\delta(s))\, ds\\
&\quad -\exp(4 w_M(t)) \int_{t_0}^t \exp(2 w_M(s))\, ds = R_M(t),
\end{aligned}
\label{inteq:6}
\end{equation}
where 
\begin{equation}
R_M(t) =   \omega^2q(t)-w_M''(t) + (w_M'(t))^2 - \exp(4 w_M(t)) \int_{t_0}^t \exp(2 w_M(s))\, ds.
\label{inteq:7}
\end{equation}
From the definitions of $w$ and $w_M$,  we see that $R_M$ is, in fact, equal to (\ref{analysis:residual}),
and therefore  is on the order of $\omega^{-2M}$.
Moreover,  because the derivatives of the $u_j$ are bounded independent of $\omega$, the same
is true of the derivatives of $R_M$.  In particular, 
\begin{equation}
\left\| R_M \right\|_\infty = \mathcal{O}\left(\frac{1}{\omega^{2M}}\right)
\quad\mbox{and}\quad
\left\| R_M' \right\|_\infty = \mathcal{O}\left(\frac{1}{\omega^{2M}}\right)
\ \ \mbox{as}\ \omega\to\infty.
\label{analysis:RMestimate}
\end{equation}
%

Expanding  the exponentials appearing in the second line of  (\ref{inteq:6}) 
in power series and rearranging the resulting expression gives us 
\begin{equation}
\begin{aligned}
\delta''(t) - &2w_M'(t)\delta'(t) + 4 \delta(t)\, \exp(4w_M(t)) \int_{t_0}^t \exp(2w_M(s))\, ds\\
&+2 \exp(4w_M(t)) \int_{t_0}^t \exp(2w_M(s))  \delta(s) \, ds=
R_M(t) +  F(t,\delta(t),\delta'(t)),
\end{aligned}
\label{inteq:9}
\end{equation}
where 
\begin{equation}
\begin{aligned}
&F(t,\delta(t),\delta'(t) ) = (\delta'(t))^2
- 8 \delta(t) \exp(4w_M(t)) \int_{t_0}^t \exp(2w_M(s))  \delta(s) \, ds\\
&\quad- \exp(4w_M(t)) \left(\frac{(4\delta(t))^2}{2} + \frac{(4\delta(t))^3}{3!} + \cdots\right)
\int_{t_0}^t \exp(2w_M(s))\exp(2\delta(s))\, ds\\
&\quad - \exp(4w_M(t)) \exp(4\delta'(t)) \int_{t_0}^t \exp(2w_M(s))\left(\frac{(2\delta(s))^2}{2} + \frac{(2\delta(s))^3}{3!}+\cdots\right)\, ds.
\end{aligned}
\label{analysis:F}
\end{equation}
Using the definition of $w_M$, we see that (\ref{inteq:9}) is equivalent to
\begin{equation}
\delta''(t) - \frac{\gamma_M''(t)}{\gamma_M'(t)}\delta'(t) + 4 (\gamma_M'(t))^2 \gamma_M(t) \delta(t) 
+ 2 (\gamma_M'(t))^2 \int_{t_0}^t \gamma_M'(s) \delta(s)\, ds
= f(t)
\label{inteq:10}
\end{equation}
with $f(t)$ taken to be  $R_M(t) +  F(t,\delta(t),\delta'(t))$.

%

We now observe that if $z$ solves the equation
\begin{equation}
\left(\gamma_M'(t)\right)^2\left(z'''(\gamma_M(t))  + 4 \gamma_M(t) z'(\gamma_M(t))  + 2 \left( z(\gamma_M(t)) - z(\gamma_M(t_0))\right)\right)
  = f(t), \ \ \ a < t < b,
\label{inteq:10.75}
\end{equation}
and we let $\delta(t) = z'(\gamma_M(t))$, then $\delta$ satisfies (\ref{inteq:10}).
Since $\gamma_M'$ is nonvanishing on $[a,b]$ for sufficiently large $\omega$ and  $\gamma_M(t_0) =0 $, 
we can divide by $(\gamma_M'(t))^2$ and introduce the new variable $u = \gamma_M(t)$ in (\ref{inteq:10.75}) to obtain
\begin{equation}
z'''(u) + 4  u z'(u)  + 2 \left(z(u)-z(0)\right)= 
f\left(\gamma_M^{-1}(u)\right)\left(\frac{1}{\gamma_M'(\gamma_M^{-1}(u))}\right)^2,
\ \ \ \gamma_M(a) < u < \gamma_M(b).
\label{inteq:12}
\end{equation}
Using the Green's function $G(t,s)$ defined in Section~\ref{section:preliminaries:green},
we can express a solution of (\ref{inteq:12}) via the formula
\begin{equation}
z(u) = \int_{\gamma_M(a)}^{\gamma_M(b)} G(u,v) f\left(\gamma_M^{-1}(v)\right)\left(\frac{1}{\gamma_M'(\gamma_M^{-1}(v))}\right)^2\, dv.
\label{inteq:12.5}
\end{equation}
We note that $G$ was chosen so that the function defined in (\ref{inteq:12.5}) satisfies the condition $z(0)=0$.
By differentiating (\ref{inteq:12.5})
%
and letting $u = \gamma_M(t)$, $v = \gamma_M(s)$ in the resulting expression, we obtain
the integral equation
\begin{equation}
\delta(t) = z'(\gamma_M(t)) = \int_a^b K(t,s)\, \left( R_M(s) + F(s,\delta(s),\delta'(s))\right)\ ds,
\label{analysis:inteq}
\end{equation}
where
\begin{equation}
K(t,s) = \frac{\partial_1 G(\gamma_M(t),\gamma_M(s))}{\gamma_M'(s)},
\label{inteq:13}
\end{equation}
$R_M$ is the residual (\ref{analysis:residual})
and $F$ is given by (\ref{analysis:F}).

We have arrived at the desired integral equation formulation (\ref{analysis:inteq}) for the function $\delta$.
Our particular choice of the Green's function $G(t,s)$ ensures
that $K(t,s)$ and its first derivatives are bounded independent of $\omega$.
In fact, from (\ref{preliminaries:airyasym1}) and  (\ref{preliminaries:airyasym2})
and the definition (\ref{analysis:gammaM}) of $\gamma_M$
it can be seen that there exists a constant $C_1$ such that 
\begin{equation}
\begin{aligned}
&\left|K(t,s)\right| \leq \frac{C_1}{\omega}\quad\mbox{and}\quad
\left|\partial_1 K(t,s)\right| \leq C_1
\end{aligned}
\label{analysis:K}
\end{equation}
for almost all $t$ and $s$.  Moreover, if we let $\widetilde{K}$ be the antiderivative
\begin{equation}
\widetilde{K}(t,s) = \int_0^t K(t,s)\, ds
\label{analysis:Ktildedef}
\end{equation}
of $K$, then we can adjust $C_1$ such that the bounds
\begin{align}
&\left|\widetilde{K}(t,s)\right| \leq \frac{C_1}{\omega^2},\quad
\left|\partial_1 \widetilde{K}(t,s)\right| \leq \frac{C_1}{\omega}
\quad\mbox{and}\quad\left|\partial_2 \widetilde{K}(t,s)\right| \leq \frac{C_1}{\omega}
\label{analysis:Ktilde}
\end{align}
also hold almost everywhere.

\end{subsection}

%
%
\begin{subsection}{Contraction mapping argument}
\label{section:analysis:contraction}

We are now in a position to give a  contraction mapping argument which shows the existence 
of a solution $\delta$ of the integral equation (\ref{analysis:inteq}) 
whose magnitude is on the order of $\omega^{-2(M+1)}$.
More explicitly, we will prove that there exists an $r>0$ such that for all sufficiently large $\omega$,
the nonlinear integral operator
\begin{equation}
T\left[\delta\right](t) = \int_a^b K(t,s)\, \left( R_M(s) + F(s,\delta(s),\delta'(s))\right)\, ds
\label{contract:1}
\end{equation}
is a contraction on a closed ball of radius $r \omega^{-2M}$ in the Banach space $C^1\left(\left[a,b\right]\right)$
endowed with the norm $\|\cdot\|_\omega$ defined in (\ref{preliminaries:wnorm}).  It will follow immediately that, assuming
$\omega$ is large enough, there exists a solution $\delta$ of (\ref{analysis:inteq}) with
$\|\delta\|_\infty \leq \frac{r}{\omega^{2(M+1)}}$.

Owing to (\ref{analysis:RMestimate}) and the definition of $\gamma_M$, we can choose constants
$C_2$ and $\omega_0$ such that 
\begin{equation}
\|\gamma_M\|_\infty \leq C_2 \omega^{\frac{2}{3}},\ \ \ \|\left(\gamma_M'\right)^2\|_\infty \leq C_2 \omega^{\frac{4}{3}},\ \ \ 
\left\| R_M\right\|_\infty \leq \frac{C_2}{\omega^{2M}}\ \ \mbox{and}\ \ 
\left\| R_M'\right\|_\infty \leq \frac{C_2}{\omega^{2M}}
\label{contract:4}
\end{equation}
whenever $\omega \geq \omega_0$.  Moreover, we note that if  $\|\delta\|_\omega \leq r \omega^{-2M}$, then 
\begin{equation}
\|\delta\|_\infty \leq \frac{r}{\omega^{2M+2}}
\ \ \ \mbox{and} \ \ \
\|\delta'\|_\infty \leq \frac{r}{\omega^{2M+1}}.
\label{contract:5}
\end{equation}
Now we use integration by parts to see that
\begin{equation}
\int_a^b K(t,s)\, R_M(s)\, ds
=
\widetilde{K}(t,b) R_M(b) - \widetilde{K}(t,a) R_M(a) -  \int_a^b \widetilde{K}(t,s)\, R_M'(s)\, ds,
\label{contract:6}
\end{equation}
where $\widetilde{K}$ is the antiderivative of $K$ defined in (\ref{analysis:Ktildedef}).
By combining (\ref{analysis:Ktilde}), (\ref{contract:4}) and (\ref{contract:6}) we obtain
the inequality 
\begin{equation}
 \left| \int_a^b K(t,s)\, R_M(s)\, ds\right| \leq  \frac{C_1 C_2 (b-a+2)}{\omega^{2M+2}},
\label{contract:7}
\end{equation}
 which holds for almost  all $t \in [a,b]$ whenever  $\omega \geq \omega_0$.
%
Now we observe that (\ref{contract:4}) and (\ref{contract:5}) imply that
\begin{equation}
\sup_{a \leq t \leq b} \left| \left(\delta'(t)\right)^2\right| \leq \frac{r^2}{\omega^{4M+2}}
\label{contract:8}
\end{equation}
and
\begin{equation}
\begin{aligned}
\sup_{a \leq t \leq b} \left|8\delta(t) \exp(4w_M(t)) \int_{t_0}^t \exp(2w_M(s))  \delta(s) \, ds \right|
&\leq  \frac{8 (b-a) C_2^2 r^2}{\omega^{4M+2}}
\end{aligned}
\label{contract:9}
\end{equation}
hold whenever $\|\delta\|_\omega \leq r \omega^{-2M}$ and $\omega \geq \omega_0$. Likewise,  we have
\begin{align}
\sup_{a \leq t \leq b}
&\left|\exp(4w_M(t)) \left(\frac{(4\delta(t))^2}{2} + \frac{(4\delta(t))^3}{3!} + \cdots\right)
\int_{t_0}^t \exp(2w_M(s))\exp(2\delta(s))\, ds \right| \nonumber \\ 
&\leq  (b-a) C_2^2 \omega^{2}  \exp\left(\frac{2r}{\omega^{2M+2}}\right)
\left( \frac{1}{2} \left(\frac{4r}{\omega^{2M+2}}\right)^2 + \frac{1}{3!} \left(\frac{4r}{\omega^{2M+2}}\right)^3
+ \cdots \right) \nonumber \\
&\leq
  (b-a) C_2^2 \omega^{2}  \exp\left(\frac{2r}{\omega^{2M+2}}\right)
\frac{16r^2}{ \omega^{4M+4}}
\left( \frac{1}{2} + \frac{1}{3!} \frac{4r}{\omega^{2M+2}}  +   \cdots \right)\label{contract:10}\\
&\leq
  (b-a) C_2^2 \omega^{2}  \exp\left(\frac{2r}{\omega^{2M+2}}\right)
\frac{16r^2}{ \omega^{4M+4}} \exp\left(\frac{4r}{\omega^{2M+2}}\right) \nonumber \\
&=
   \frac{ 16 (b-a)   C_2^2r^2}{ \omega^{4M+2}}  \exp\left(\frac{6r}{\omega^{2M+2}}\right)\nonumber 
\end{align}
and
\begin{align}
\sup_{a\leq t\leq b}&\left|\exp(4w_M(t)) \exp(4\delta(t))\int_{t_0}^t \exp(2w_M(s)) \left( \frac{(2\delta(s))^2}{2} + \frac{(2\delta(s))^3}{3!}    + \cdots \right) \, ds\right|\nonumber \\
&\leq 
  \omega^2 (b-a) C_2^2\exp\left(\frac{4r}{\omega^{2M+2}}\right) \left( 
\frac{1}{2} \left(\frac{2r}{\omega^{2M+2}} \right)^2 + 
\frac{1}{3!} \left(\frac{2r}{\omega^{2M+2}} \right)^3
+\cdots
\right)\label{contract:11}\\
&\leq
 \frac{4  (b-a)C_2^2r^2}{\omega^{4M+2}}  \exp\left(\frac{6r}{\omega^{2M+2}}\right)\nonumber 
\end{align}
whenever $\omega \geq \omega_0$ and  $\|\delta\|_\omega \leq r$.   By combining  
(\ref{contract:7}) through (\ref{contract:11})
and making use of  (\ref{analysis:K}), we see that
\begin{equation}
\begin{aligned}
\left\|T\left[\delta \right]\right\|_\infty
\leq &
\frac{(b-a+2)C_1C_2}{\omega^{2M+2}} + \frac{r^2 (b-a)C_1}{\omega^{4M+3}}+ \frac{8 (b-a)^2 C_1 C_2^2 r^2}{\omega^{4M+3}}\\
 +&\frac{ 16 (b-a)^2  C_1 C_2^2 r^2}{ \omega^{4M+3}}  \exp\left(\frac{6r}{\omega^{2M+2}}\right)+
 \frac{4(b-a)^2C_1C_2^2r^2  }{\omega^{4M+3}}  \exp\left(\frac{6r}{\omega^{2M+2}}\right)
\end{aligned}
\label{contract:12}
\end{equation}
whenever $\omega \geq \omega_0$ and  $\|\delta\|_\omega \leq r$.   We can rearrange (\ref{contract:12})
as
\begin{equation}
\left\|T\left[\delta \right]\right\|_\infty
\leq \frac{1}{\omega^{2M+2}}\left((b-a+2)C_1 C_2  + D(r,\omega)
\right),
\label{contract:13}
\end{equation}
where $D(r,\omega) \to 0$ as $\omega\to\infty$ for any fixed $r$. An essentially identical argument  
gives us the bound
\begin{equation}
\left\|T\left[\delta \right]'\right\|_\infty
\leq \frac{1}{\omega^{2M+1}}\left( (b-a+2)C_1 C_2 + D(r,\omega) \right)
\label{contract:14}
\end{equation}
on the derivative $T\left[\delta\right]'$ of $T\left[\delta\right]$.  
We now choose $r=2 (b-a+2)C_1 C_2$.  Then, for any sufficiently large $\omega$, 
$T$ preserves the ball of radius $r \omega^{-2M}$ centered at $0$ in the space $C^1\left(\left[a,b\right]\right)$
endowed with the norm (\ref{preliminaries:wnorm}).

It remains to show that $T$ is a contraction for sufficiently large $\omega$.   To that end,
we observe that the Fr\'echet derivative of $T$ at the point $\delta$ is the linear operator
\begin{equation}
D_\delta T\left[h\right](t) =  \int_a^b K(t,s)\, H(s,\delta(s), \delta'(s), h(s), h'(s))\, ds,
\label{contract:16}
\end{equation}
where $H(s,\delta(s),\delta'(s), h(s), h'(s))$ is 
\begin{align}
&2 \delta'(t) h'(t)- \left(8 \exp(4w_M(t))  \int_{t_0}^t \exp(2w_M(s)) \delta(s)\, ds\right) h(t) \nonumber\\
&- 8 \delta(t) \exp(4w_M(t)) \int_{t_0}^t \exp(2w_M(s)) h(s)\, ds\nonumber\\
&- 4\left(\exp(4w_M(t)) \left( 4\delta(t) + \frac{(4\delta(t))^2}{2} + \cdots \right) 
\int_{t_0}^t \exp(2w_M(s)) \exp(2\delta(s)) \, ds \right)h(t)\\
&- 2\exp(4w_M(t)) \left( \frac{(4\delta(t))^2}{2} + \frac{(4\delta(t))^3}{3!}+ \cdots \right) \int_{t_0}^t \exp(2w_m(s)) \exp(2\delta(s)) h(s)\, ds\nonumber\\
&- 4\left(\exp(4w_M(t)) \exp(4\delta(t)) \int_{t_0}^t \exp(2w_M(s)) \left(\frac{(2\delta(s))^2}{2} + \frac{(2\delta(s))^3}{3!} + \cdots\right)\, ds
\right)h(t)\nonumber\\
&- 2\exp(4w_M(t)) \exp(4\delta(t)) \int_{t_0}^t \exp(2w_M(s)) \left(2\delta(s) + \frac{(2\delta(s))^2}{2!} + \cdots\right)h(s)\, ds.\nonumber
\label{contract:17}
\end{align}
Using (\ref{contract:4}) and (\ref{contract:5}), we see that
%
\begin{align}
\|H(s,\,&\delta(s),\,\delta'(s),\,h(s),\,h'(s))\|_\infty
\leq 
\frac{2r}{\omega^{2M+1}} + \frac{8(b-a)C_2^2 r}{\omega^{2M+2}}
+ \frac{8(b-a)C_2^2 r}{\omega^{2M+2}} \nonumber\\
&+
\frac{16 (b-a)C_2^2 r}{\omega^{2M+2}} \exp\left(\frac{6r}{\omega^{2M+2}}\right)
+
\frac{16(b-a)C_2^2 r^2}{\omega^{4M+4}} \exp\left(\frac{6r}{\omega^{2M+2}}\right)+\\
&+
\frac{8(b-a)C_2^2r^2 }{\omega^{4M+4}}\exp\left(\frac{6r}{\omega^{2M+2}}\right)
+\frac{4(b-a)C_2^2r }{\omega^{2M+2}}\exp\left(\frac{6r}{\omega^{2M+2}}\right)\nonumber
\label{contract:18}
\end{align}
%
whenever $\left\|\delta\right\|_\omega \leq r$, $\omega \geq \omega_0$ and $\|h\|_\omega=1$.
In particular, there exists a function $E(\omega)$ which is bounded
for all $\omega \geq \omega_0$ and such that
\begin{equation}
\|H(s,\delta(s),\delta'(s),h(s),h'(s))\|_\infty \leq \frac{ E(\omega)}{\omega^{2M+1}}
\label{contract:19}
\end{equation}
when $\left\|\delta\right\|_\omega \leq r$, $\omega \geq \omega_0$ and $\|h\|_\omega=1$.
From (\ref{analysis:K}) and (\ref{contract:19}), we have that
\begin{equation}
\left\|D_\delta T\left[h\right] \right\|_\infty 
\leq \frac{ C_1 (b-a)E(\omega)}{\omega^{2M+2}}\ \ \ \mbox{and} \ \ \ 
\left\|D_\delta T\left[h\right]' \right\|_\infty 
\leq \frac{ C_1 (b-a)E(\omega)}{\omega^{2M+1}}.
\end{equation}
In particular, the $\|\cdot\|_\omega$ operator norm of the Fr\'echet derivative $D_\delta T$ of $T$ at the
point $\delta$ is less than $1$ whenever $\delta$ is in the ball of radius $r$ centered at $0$
in $C^1\left(\left[a,b\right]\right)$ and $\omega$ is sufficiently large.  It follows
that $T$ is a contraction and our proof of the existence of a slowly-varying
solution of the Airy-Kummer equation is complete.

\end{subsection}

\end{section}

\begin{section}{Convergence of the Newton-Kantorovich Method}
\label{section:newton}

In this section, we discretize the Airy-Kummer equation (\ref{introduction:airykummer})
 over an interval of the form $[-a_0,a_0]$ via a Chebyshev  spectral method 
and apply  Newton-Kantorovich method to the resulting system of nonlinear algebraic equations.
We show that if the first order approximate $\gamma_0$ is used to form an initial guess, then 
the Newton iterates converge to a vector representing the slowly-varying solution of 
the Airy-Kummer equation  whose existence was established in the preceding section provided 
$\omega$ is sufficiently large and $a_0$ is sufficiently small.
In the course of our proof, we will see that  the number of nodes $k$ used to discretize the 
Airy-Kummer equation must be relatively small (we take $k=16$ in the experiments of this paper,
and this causes no difficulties).  However, because the Chebyshev spectral
discretizations of the various functions which arise in the course of the proof converge as $a_0$ goes to $0$, 
this does not  impose a significant limitation on the applicability of our algorithm.

We divide our argument into four parts.  In Subsection~\ref{section:newton:spectral}, we 
form the spectral discretization of the Airy-Kummer equation. 
We then show in Subsection~\ref{section:newton:frechet} 
that the Fr\'echet derivative of the operator in question  is invertible 
provided $a_0$ is sufficiently small and $\omega$ sufficiently large.
In Subsection~\ref{section:newton:lipschitz}, we give a Lipschitz bound on the Fr\'echet derivative.
Finally, we conclude our proof by applying  the Newton-Kantorovich theorem
 in Subsection~\ref{section:newton:nk}.

\begin{subsection}{Spectral discretization of the Airy-Kummer equation}
\label{section:newton:spectral}

We form a spectral discretization of the Airy-Kummer equation (\ref{introduction:airykummer}) by representing the unknown
solution $\gamma$  and the coefficient $q$ via the vectors
\begin{equation}
\boldsymbol{\gamma} = \left(\begin{array}{cccc}
\gamma(t_0) & \gamma(t_1) & \cdots & \gamma(t_k)
\end{array}\right)^\top
\end{equation}
and
\begin{equation}
\boldsymbol{q} = \left(\begin{array}{ccccc}
q(t_1) & q(t_2) & \cdots & q(t_k)
\end{array}\right)^\top
\label{newton:q}
\end{equation}
of their values at the extremal Chebyshev nodes $t_1,\ldots,t_k$ on the interval $[-a_0,a_0]$, 
replacing the differential operators with Chebyshev spectral differentiation matrices
and requiring that the resulting system of semidiscrete equations hold at the nodes $t_1,\ldots,t_k$.
This procedure yields the system of nonlinear algebraic equations $R\left(\boldsymbol{\gamma}\right)=0$,
where  $R$ is the mapping $\mathbb{R}^k \to \mathbb{R}^k$ given by the matrix
\begin{equation} 
\label{newton:R}
R\left(\boldsymbol{\gamma}\right) = 
\omega^2 \boldsymbol{q} 
        - \boldsymbol{\gamma} \circ 
        \left(
            \frac{\mathscr{D}_k}{a_0} \boldsymbol{\gamma}
        \right)^{\circ 2}
        + \frac{3}{4} 
        \left(
            \left(
                \frac{\mathscr{D}_k}{a_0}
            \right)^2 \boldsymbol{\gamma}
        \right)^{\circ 2} 
        \circ 
        \left(
            \frac{\mathscr{D}_k}{a_0} \boldsymbol{\gamma}
        \right)^{\circ -2}
        - \frac{1}{2} 
        \left(
            \left( \frac{\mathscr{D}_k}{a_0} \right)^3 \boldsymbol{\gamma} 
        \right)
        \circ 
        \left(
            \frac{\mathscr{D}_k}{a_0} \boldsymbol{\gamma} 
        \right)^{\circ -1}.
\end{equation}
%
The Fr\'echet derivative
of $R$ at  $\boldsymbol{x}$ is the linear mapping $\mathbb{R}^k \to \mathbb{R}^k$ 
given by the matrix
\begin{equation} 
\label{newton:DR}
 \begin{split}
  D_{\boldsymbol{x}} R = 
  &- \left(\diag \left( \frac{\mathscr{D}_k}{a_0} \boldsymbol{x} \right)
  \right)^2 - 2\diag \left(\boldsymbol{x}\right) \diag\left(\frac{\mathscr{D}_k}{a_0}\boldsymbol{x}\right) \frac{\mathscr{D}_k}{a_0}\\
  & 
  -\frac{3}{2} 
  \left( 
      \diag \left( \left( \frac{\mathscr{D}_k}{a_0} \right)^2 \boldsymbol{x} \right)
  \right)^2 
  \left(
      \diag \left( \frac{\mathscr{D}_k}{a_0} \boldsymbol{x} \right)
  \right)^{-2} \frac{\mathscr{D}_k}{a_0}\\
  &
  +\frac{1}{2}
      \diag \left( \left( \frac{\mathscr{D}_k}{a_0} \right)^3 \boldsymbol{x} \right)
  \left(
      \diag \left( \frac{\mathscr{D}_k}{a_0} \boldsymbol{x} \right)
  \right)^{-2} \frac{\mathscr{D}_k}{a_0}\\
        &+\frac{3}{2} \diag
        \left( 
            \left( \frac{\mathscr{D}_k}{a_0} \right)^2 \boldsymbol{x} 
        \right)  
        \left(
            \diag\left( \frac{\mathscr{D}_k}{a_0} \boldsymbol{x} \right)
        \right)^{-2} 
        \left( \frac{\mathscr{D}_k}{a_0} \right)^2 \\
        &- \frac{1}{2} 
        \left(
            \diag\left( \frac{\mathscr{D}_k}{a_0} \boldsymbol{x} \right)
        \right)^{-1} 
        \left( \frac{\mathscr{D}_k}{a_0} \right)^3.
    \end{split}
\end{equation}

We let $\boldsymbol{\gamma_0}$ be the spectral discretization
\begin{equation}
\boldsymbol{\gamma_0} = \left(\begin{array}{ccccc}
\gamma_0(t_1) & \gamma_0(t_2) & \cdots & \gamma_0(t_k)
\end{array}\right)^\top
\label{newton:gamma0}
\end{equation}
of the first order asymptotic approximation $\gamma_0$ defined via (\ref{introduction:gamma0}).  Owing to our assumptions
on $q$, there exists a smooth function $\tau$ such that $\tau(t) = \mathcal{O}\left(t^2\right)$ as $t \to 0$ 
and  $\gamma_0(t) = \omega^{\frac{2}{3}} \left(t + \tau(t)\right)$.  We now write
\begin{equation}
\boldsymbol{\gamma_0} = \omega^{\frac{2}{3}} \left( a_0 \boldsymbol{t} + \boldsymbol{\tau} + \boldsymbol{\epsilon_1}\right)
\quad\mbox{and}\quad
\frac{\mathscr{D}_k}{a_0} \boldsymbol{\gamma_0} = \omega^{\frac{2}{3}} \left( \boldsymbol{1} + \boldsymbol{\tau'} + \boldsymbol{\epsilon_2}\right),
\label{newton:tau}
\end{equation}
where:
\renewcommand\labelitemi{\tiny$\bullet$}
\begin{itemize}
\item
%
$\begin{aligned}\boldsymbol{t} = \left(\begin{array}{ccccc}
t_1^{\mbox{\tiny cheb}} & t_2^{\mbox{\tiny cheb}} & \cdots & t_k^{\mbox{\tiny cheb}}
\end{array}\right)^\top\ \mbox{and}\ 
\boldsymbol{1} = \left(\begin{array}{ccccc}
1 & 1 & \cdots & 1
\end{array}\right)^\top\end{aligned}$;
%
\item $\boldsymbol{\tau}$ and $\boldsymbol{\tau'}$ are the spectral discretizations
\begin{equation}
\boldsymbol{\tau} = \left(\begin{array}{cccc} \tau(t_1) & \tau(t_2) & \cdots & \tau(t_k)\end{array}\right)^\top\quad\mbox{and}\quad
\boldsymbol{\tau'} = \left(\begin{array}{cccc} \tau'(t_1) & \tau'(t_2) & \cdots & \tau'(t_k)\end{array}\right)^\top
\end{equation}
of the function $\tau$ and $\tau'$; and 
\item $\boldsymbol{\epsilon_1}$ and $\boldsymbol{\epsilon_2}$  are vectors which account for discretization error.
\end{itemize}
Since $\tau(t) = \mathcal{O}\left(t^2\right)$,  we have
\begin{equation}
\left\|\boldsymbol{\tau}\right\|_\infty = \mathcal{O}\left(a_0^2\right)\quad\mbox{and}\quad
\left\|\boldsymbol{\tau'}\right\|_\infty = \mathcal{O}\left(a_0\right)
 \ \mbox{as}\ \ a_0\to 0,
\label{newton:tauasym}
\end{equation}
while the infinite differentiability of $\gamma_0$ implies
\begin{equation}
\left\|\boldsymbol{\epsilon_1}\right\|_\infty = \mathcal{O}\left(a_0^{k}\right)
\quad\mbox{and}\quad
\left\|\boldsymbol{\epsilon_2}\right\|_\infty = \mathcal{O}\left(a_0^{k-1}\right)
\ \mbox{as} \ \ a_0 \to 0.
\label{newton:disclimit}
\end{equation}
Moreover, because of the form of $\gamma_0$ and our assumption that the derivatives of $q$ are bounded independent
of $\omega$, the asymptotic estimates (\ref{newton:tauasym}) and (\ref{newton:disclimit})
hold uniformly in $\omega$.
We close this subsection by noting that, because the residual $R_0$ defined in (\ref{analysis:residual}) is bounded independent of $\omega$
and $R\left(\boldsymbol{\gamma_0}\right)$ discretizes this quantity,
we will have
\begin{equation}
R\left(\boldsymbol{\gamma_0}\right) = \mathcal{O}\left(1\right)
\ \mbox{as}\ \omega\to\infty,
\label{newton:residual}
\end{equation}
for any fixed $a_0$ which is sufficiently small so that the spectral discretization
of $R\left(\boldsymbol{\gamma_0}\right)$  is accurate.

\end{subsection}

\begin{subsection}{The Fr\'echet derivative of $R$ at $\boldsymbol{\gamma_0}$}
\label{section:newton:frechet}

It follows from (\ref{newton:DR}) and (\ref{newton:tau}) that the matrix representing
the Fr\'echet derivative of $R$ at  $\boldsymbol{\gamma_0}$ can be written
in the form
\begin{equation}
\begin{aligned}
D_{\boldsymbol{\gamma_0}} R
= 
-\omega^{\frac{4}{3}} \left(
\mathscr{I}_k + 2 \diag(\boldsymbol{t}) \mathscr{D}_k+S
\right)+ T,
\end{aligned}
\label{newton:op}
\end{equation}
where 
\begin{equation}
\begin{aligned}
S = 
&\left(\diag(\boldsymbol{\tau'})\right)^2  + 
\left(\diag(\boldsymbol{\epsilon_2})\right)^2 +
2 \diag(\boldsymbol{\tau'})\diag(\boldsymbol{\epsilon_2}) +  
2 \diag(\boldsymbol{\tau'}) + 2\diag(\boldsymbol{\epsilon_2}) + \\
& 2\diag(\boldsymbol{t}) \diag(\boldsymbol{\tau'})\mathscr{D}_k +
2\diag(\boldsymbol{t}) \diag(\boldsymbol{\epsilon_2}) \mathscr{D}_k+\\
&2 \diag(\boldsymbol{\tau})\frac{\mathscr{D}_k}{a_0} +
2 \diag(\boldsymbol{\tau})\diag(\boldsymbol{\tau'})\frac{\mathscr{D}_k}{a_0} + 
2 \diag(\boldsymbol{\tau})\diag(\boldsymbol{\epsilon_2})\frac{\mathscr{D}_k}{a_0}+\\
& 2 \diag(\boldsymbol{\epsilon_1})\frac{\mathscr{D}_k}{a_0} + 2 \diag(\boldsymbol{\epsilon_1})\diag(\boldsymbol{\tau'})\frac{\mathscr{D}_k}{a_0} + 
2 \diag(\boldsymbol{\epsilon_1})\diag(\boldsymbol{\epsilon_2})\frac{\mathscr{D}_k}{a_0}
\end{aligned}
\label{newton:S}
\end{equation}
and the operator norm of the matrix $T$ is bounded independent of $\omega$.
The  matrix $\mathscr{I}_k + 2\diag(\boldsymbol{t})\mathscr{D}_k$ is invertible and
well-conditioned  as long as $k$ is of moderate size.   From (\ref{newton:tauasym}) and (\ref{newton:disclimit}),
it is clear that the operator norm of $S$ goes to $0$ as $a_0 \to 0$.  
It follows from this and a standard Neumann series argument that the operator
\begin{equation}
\mathscr{I}_k + 2 \diag(\boldsymbol{t}) \mathscr{D}_k+S
\label{newton:invert}
\end{equation}
is invertible for all sufficiently small $a_0$. 
Since the operator norm of $S$ is bounded independent of $\omega$,
a second routine Neumann series argument implies that $D_{\boldsymbol{\gamma_0}} R$
is invertible provided (\ref{newton:invert}) is invertible and $\omega$ is sufficiently large.
Moreover, it is clear from (\ref{newton:op}) that 
\begin{equation}
\left\|\left(D_{\boldsymbol{\gamma_0}} R\right)^{-1}\right\|_\infty = \mathcal{O}\left(\omega^{-\frac{4}{3}}\right)
\ \mbox{as}\ \omega\to\infty.
\label{newton:frechetbound}
\end{equation}

\end{subsection}

\begin{subsection}{A Lipschitz bound for $\left\|D_{\boldsymbol{x}}R- D_{\boldsymbol{y}}R\right\|_\infty$}
\label{section:newton:lipschitz}

From (\ref{newton:DR}), we see that the Fr\'echet derivative of $R$ at $\boldsymbol{v}$
can be written as 
\begin{equation} 
\label{newton:dr2}
D_{\boldsymbol{v}} R  = 
    \sum_{j=0}^{3} \diag\left( E_j\left(\boldsymbol{v}\right) \right) 
    \left( \frac{\mathscr{D}_k}{a_0} \right)^j,
\end{equation}
where the $E_j$ are the nonlinear mappings $\mathbb{R}^k \to \mathbb{R}^k$ defined
via the following formulas:
\begin{equation}
    \begin{split}
        E_0(\boldsymbol{v})
        &= 
        - \left( 
            \frac{\mathscr{D}_k}{a_0} \boldsymbol{v} 
        \right)^{\circ 2}, \\
        E_1(\boldsymbol{v})
        &= 
        -2 \boldsymbol{v} \circ \left( \frac{\mathscr{D}_k}{a_0} \boldsymbol{v} \right)
        -\frac{3}{2}
        \left( 
            \left( \frac{\mathscr{D}_k}{a_0} \right)^2 \boldsymbol{v} 
        \right)^{\circ 2} 
        \circ 
        \left( \frac{\mathscr{D}_k}{a_0} \boldsymbol{v} \right)^{\circ -3} 
        + \frac{1}{2}
        \left( 
            \left( \frac{\mathscr{D}_k}{a_0} \right)^3 \boldsymbol{v} 
        \right) 
        \circ 
        \left( 
            \frac{\mathscr{D}_k}{a_0} \boldsymbol{v}
        \right)^{\circ -2}, \\
        E_2(\boldsymbol{v})
        &= 
        \frac{3}{2} 
        \left( 
            \left( \frac{\mathscr{D}_k}{a_0} \right)^2 \boldsymbol{v} 
        \right) 
        \circ \left( \frac{\mathscr{D}_k}{a_0} \boldsymbol{v} \right)^{\circ -2} \quad\mbox{and}
        \\
        E_3(\boldsymbol{v})
        &= 
        - \frac{1}{2} \left( \frac{\mathscr{D}_k}{a_0} \boldsymbol{v} \right)^{\circ -1}.
    \end{split}
\end{equation}
Given any $\eta >0$, it follows from (\ref{newton:dr2}) that
\begin{equation} \
\label{newton:xminusy}
    \begin{split}
\left\|
 D_{\boldsymbol{x}} R-  D_{\boldsymbol{y}} R
        \right\|_\infty
        &\le
        \sum_{j=0}^3 
        \left(
            \left\|
                E_j\left(\boldsymbol{x}\right) - E_j\left(\boldsymbol{y}\right)
            \right\|_\infty
            \left\| \frac{\mathscr{D}_k}{a_0} \right\|_\infty^j
        \right) \\
        &\le 
        \left(
            \sum_{j=0}^3 
                \left\| \frac{\mathscr{D}_k}{a_0} \right\rVert_{\infty}^j \
                \sup_{\boldsymbol{v} \in B_{\eta}(\boldsymbol{\gamma_0})} 
                \left\| D_{\boldsymbol{v}} E_j \right\|_\infty
            \right) 
            \left\| \boldsymbol{x} - \boldsymbol{y} \right\|_\infty
    \end{split}
\end{equation}
whenever  $\boldsymbol{x}, \boldsymbol{y} \in B_\eta\left(\boldsymbol{\gamma_0}\right)$.
Simple calculations shows that the Fr\'echet derivatives of the $E_j$ are given by the matrices
%
\begin{align}
\label{newton:DEs}
        D_{\boldsymbol{v}} E_0 
        &= 
        -2 \ \diag 
        \left( \frac{\mathscr{D}_k}{a_0} \boldsymbol{v} \right) 
        \frac{\mathscr{D}_k}{a_0}, \nonumber \\
        D_{\boldsymbol{v}} E_1
        &= 
        - 2 \ 
        \left(
            \diag \left( \frac{\mathscr{D}_k}{a_0} \boldsymbol{v} \right) 
            + \diag (\boldsymbol{v}) \frac{\mathscr{D}_k}{a_0} 
        \right)  \nonumber  \\
        &\quad 
        - \frac{3}{2} 
        \Bigg(
            -3 \ 
            \left(
                \diag \left( \left( \frac{\mathscr{D}_k}{a_0} \right)^2 \boldsymbol{v} \right)
            \right)^2  
            \left(
                \diag \left( \frac{\mathscr{D}_k}{a_0} \boldsymbol{v} \right) 
            \right)^{-4} 
            \frac{\mathscr{D}_k}{a_0}  \nonumber \\
            &\qquad\qquad
            + 2 \ \diag 
            \left( 
                \left( \frac{\mathscr{D}_k}{a_0} \right)^2 \boldsymbol{v}
            \right)  
            \left(
                \diag \left( \frac{\mathscr{D}_k}{a_0} \boldsymbol{v} \right)
            \right)^{-3} 
            \left( \frac{\mathscr{D}_k}{a_0} \right)^2 
        \Bigg) \\
        &\quad
        + \frac{1}{2} 
        \left( 
            -2 \ \diag
            \left( 
                \left( \frac{\mathscr{D}_k}{a_0} \right)^3 \boldsymbol{v}
            \right) 
            \left( 
                \diag \left( \frac{\mathscr{D}_k}{a_0} \boldsymbol{v} \right)
            \right)^{-3} 
            \frac{\mathscr{D}_k}{a_0}
            + 
            \left(
                \diag \left( \frac{\mathscr{D}_k}{a_0} \boldsymbol{v} \right) 
            \right)^{-2} 
            \left( \frac{\mathscr{D}_k}{a_0} \right)^3 
        \right),  \nonumber \\
        D_{\boldsymbol{v}} E_2
        &=
        \frac{3}{2} 
        \left(
            -2 \ \diag 
            \left( 
                \left( \frac{\mathscr{D}_k}{a_0} \right)^2 \boldsymbol{v}
            \right)
            \left( 
                \diag \left( \frac{\mathscr{D}_k}{a_0} \boldsymbol{v} \right)
            \right)^{-3} 
            \frac{\mathscr{D}_k}{a_0}
            + 
            \left( 
                \diag \left( \frac{\mathscr{D}_k}{a_0} \boldsymbol{v} \right)
            \right)^{-2} 
            \left( \frac{\mathscr{D}_k}{a_0} \right)^2
        \right)\quad\mbox{and}  \nonumber \\
        D_{\boldsymbol{v}} E_3
        &=
        \frac{1}{2} 
        \left(
            \diag \left( \frac{\mathscr{D}_k}{a_0} \boldsymbol{v} \right)
        \right)^{-2} 
        \frac{\mathscr{D}_k}{a_0}. \nonumber 
\end{align}
%
It follows from (\ref{newton:tau}), (\ref{newton:tauasym}) and (\ref{newton:disclimit})
that for sufficiently small $a_0$, we can choose a constant   $C_3$ such that 
$\left\|\boldsymbol{\gamma_0}\right\|_\infty \leq C_3 \omega^{\frac{2}{3}}$
and each entry $v_j$ of the vector $\frac{\mathscr{D}_k}{a_0} \boldsymbol{\gamma_0}$ satisfies the inequality
\begin{equation}
\frac{C_3}{2} \omega^{\frac{2}{3}}  \leq v_j \leq C_3 \omega^{\frac{2}{3}}.
\label{newton:vj}
\end{equation}
It follows that if 
\begin{equation}
\eta < \frac{C_3 a_0}{4 \left\|\mathscr{D}_k\right\|_\infty} \omega^{\frac{2}{3}},
\label{newton:eta}
\end{equation}
then
\begin{equation}
\left\|\boldsymbol{v} \right\|_\infty\leq \left(1 + \frac{a_0}{4 \left\|\mathscr{D}_k\right\|_\infty}\right) C_3 \omega^{\frac{2}{3}},
\quad
\left\|\frac{\mathscr{D}_k}{a_0}\boldsymbol{v} \right\|_\infty\leq 
 \frac{5 C_3}{4}\omega^{\frac{2}{3}}
\quad\mbox{and}\quad
\left\|\left(\frac{\mathscr{D}_k}{a_0}\boldsymbol{v}\right)^{\circ -1} \right\|_\infty\leq 
 \frac{C_3}{4}\omega^{\frac{2}{3}}
\label{newton:vectorbounds}
\end{equation}
for all $\boldsymbol{v} \in B_\eta\left(\boldsymbol{\gamma_0}\right)$.  
Now (\ref{newton:DEs}) and (\ref{newton:vectorbounds}) imply that 
for all $\boldsymbol{v} \in B_\eta\left(\boldsymbol{\gamma_0}\right)$, 
\begin{equation}
\left\|D_{\boldsymbol{v}} E_0 \right\|_\infty = \mathcal{O}\left(\omega^{\frac{2}{3}}\right)\quad\mbox{and}\quad
\left\|D_{\boldsymbol{v}} E_1 \right\|_\infty = \mathcal{O}\left(\omega^{\frac{2}{3}}\right)
\ \mbox{as}\ \omega\to \infty
\label{newton:DE01}
\end{equation}
while
\begin{equation}
\left\|D_{\boldsymbol{v}} E_2 \right\|_\infty = \mathcal{O}\left(\omega^{-\frac{4}{3}}\right)\quad\mbox{and}\quad
\left\|D_{\boldsymbol{v}} E_3 \right\|_\infty = \mathcal{O}\left(\omega^{-\frac{4}{3}}\right)
\ \mbox{as}\ \omega\to \infty.
\label{newton:DE23}
\end{equation}
From (\ref{newton:xminusy}), (\ref{newton:DE01}) and (\ref{newton:DE23}), we see that
\begin{equation}
\left\|
 D_{\boldsymbol{x}} R-  D_{\boldsymbol{y}} R        \right\|_\infty = 
\mathcal{O}\left(\omega^{\frac{2}{3}}\right) \left\|\boldsymbol{x}-\boldsymbol{y}\right\|_\infty
\ \ \ \mbox{as}\ \omega\to\infty
\label{newton:lipschitzbound}
\end{equation}
for all $\boldsymbol{x}$ and $\boldsymbol{y}$  in $B_\eta\left(\boldsymbol{\gamma_0}\right)$
provided (\ref{newton:eta}) holds.

\end{subsection}

\begin{subsection}{Application of the Newton-Kantorovich theorem}
\label{section:newton:nk}

We are now in a position to apply the Newton-Kantorovich theorem.
We choose $a_0$ to be sufficiently small that 
(\ref{newton:residual}) holds, the operator (\ref{newton:invert}) is invertible and 
such that (\ref{newton:vj}) holds.
From the discussion in the preceding subsection, we know that if
$\eta$ satisfies (\ref{newton:eta}), then 
the  entries of the vector $\frac{\mathscr{D}_k}{a_0} \boldsymbol{v}$ are
bounded away from $0$ for all $\boldsymbol{v} \in B_\eta\left(\boldsymbol{\gamma_0}\right)$.
It follows from this and (\ref{newton:R}) that the mapping $R$ is continuously
differentiable on $B_\eta\left(\boldsymbol{\gamma_0}\right)$, and we
take the open set $\Omega$ in Theorem~\ref{theorem:nk} to be $B_\eta\left(\boldsymbol{\gamma_0}\right)$.
Condition (2) of the theorem is obviously satisfied. 
We showed in Subsection~\ref{section:newton:frechet} that the Fr\'echet derivative of $R$ at $\boldsymbol{\gamma_0}$
is invertible provided $\omega$ is sufficiently large; that is to say, condition (1) of the theorem is satisfied provided
$\omega$ is large enough.

We now choose the parameters $\eta$ and $\lambda$ as follows:
\begin{equation}
\lambda = \left\| \left(D_{\boldsymbol{\gamma_0}} R\right)^{-1} 
R\left(\boldsymbol{\gamma_0}\right)
\right\|_\infty
\quad\mbox{and}\quad
\eta = \omega^{\frac{1}{3}}.
\label{newton:etalambda}
\end{equation}
The inequality (\ref{newton:eta}) clearly holds for sufficiently large $\omega$, 
and our choice of $\lambda$ means that condition (4) of the theorem is obviously satisfied.
Now, combining (\ref{newton:residual}) and (\ref{newton:frechetbound}) yields
\begin{equation}
\lambda = \mathcal{O}\left(\omega^{-\frac{4}{3}}\right) \ \mbox{as}\ \omega\to\infty,
\label{newton:lambda_asym}
\end{equation}
and it is immediate from (\ref{newton:eta}) and (\ref{newton:lambda_asym}) 
that condition (3) of the theorem is satisfied for sufficiently large $\omega$.
It remains only to show that condition (5) holds.  Combining 
(\ref{newton:frechetbound}) and (\ref{newton:lipschitzbound}) shows that 
\begin{equation}
\left\|
\left(D_{\boldsymbol{\gamma_0}} R\right)^{-1}
\left(D_{\boldsymbol{x}} R - D_{\boldsymbol{y}} R\right)
\right\|_\infty
=  \mathcal{O}\left(\omega^{-\frac{2}{3}}\right) \left\| \boldsymbol{x} - \boldsymbol{y}\right\|_\infty
\ \mbox{as}\ \omega\to\infty
\end{equation}
whenever $\boldsymbol{x}$ and $\boldsymbol{y}$ are elements of $B_\eta\left(\boldsymbol{\gamma_0}\right)$.
It follows from this and our choice of $\eta$ that condition (5) is satisfied for sufficiently large $\omega$.

Having established that all of the requirements of Theorem~\ref{theorem:nk} are satisfied, 
it now follows that, when $\omega$ is sufficiently large, there is a unique solution
$\boldsymbol{x^*}$ of $R\left(\boldsymbol{x}\right) = 0$ in 
the ball  $B_\eta\left(\boldsymbol{\gamma_0}\right)$, and 
the sequence $\{\boldsymbol{x_j}\}$ of Newton-Kantorovich iterates 
generated by the initial guess $\boldsymbol{x_0} = \boldsymbol{\gamma_0}$ converge to $\boldsymbol{x^*}$.
Because the slowly-varying solution of the Airy-Kummer equation
whose existence was established in Section~\ref{section:analysis} converges to $\gamma_0$
as $\omega\to\infty$, the vector
$\boldsymbol{\gamma^*}$ discretizing it 
will lie in the ball $B_\eta\left(\boldsymbol{\gamma_0}\right)$ for all sufficiently large $\omega$.  
Moreover, since the equation
$R\left(\boldsymbol{x}\right)=0$ discretizes the Airy-Kummer equation which
the slowly-varying solution  satisfies, $R\left(\boldsymbol{\gamma^*}\right) \approx 0$.
It follows that $\boldsymbol{x^*}$, which is the unique solution of the discretized equation
in the ball $B_\eta\left(\boldsymbol{\gamma_0}\right)$, 
must agree closely with the discretization $\boldsymbol{\gamma^*}$ of the  slowly-varying solution.
Finally, we observe that 
(\ref{preliminaries:nkconclusion}) together with (\ref{newton:etalambda})
and (\ref{newton:lambda_asym}) give us the estimate
\begin{equation}
\left\|\boldsymbol{x_j} - \boldsymbol{\gamma^*} \right\|_\infty = \mathcal{O}\left(\omega^{\frac{1}{3} - \frac{5}{3}  2^j} \right)
\ \mbox{as}\  \omega\to\infty
\label{newton:rate}
\end{equation}
on the rate at which the Newton iterates converge to  $\boldsymbol{x^*} \approx \boldsymbol{\gamma^*}$.

\begin{remark}
Given  any $\epsilon >0$, we can  choose $\eta = \omega^{\frac{2}{3}-\epsilon}$ and the proof
given here is still valid.   In this event, 
\begin{equation}
\left\|\boldsymbol{x_j} - \boldsymbol{\gamma^*} \right\|_\infty = \mathcal{O}\left(\omega^{\left(\frac{2}{3}-\epsilon\right) - 2^{j}\left(2-\epsilon\right)} \right)
\ \mbox{as}\  \omega\to\infty.
\end{equation}
It seems likely that if constants are carefully  tracked and  accounted for,
then it will emerge that the   parameter $\eta$ can be chosen to be
on the order of $\omega^{\frac{2}{3}}$.  If  this is the case, then the  Newton iterates converge to $\boldsymbol{\gamma^*}$  at the rate
$\omega^{\frac{2}{3}-2^{j+1}}$.  
\end{remark}

\end{subsection}

\end{section}

\begin{section}{Numerical construction of Airy phase functions}
\label{section:algorithm}

In this section, we describe our numerical method for calculating the  slowly-varying Airy phase 
function $\gamma$ whose existence was established in Section~\ref{section:analysis}.
Our algorithm proceeds in three stages.  In the first stage, we numerically compute the values of the approximation $\gamma_0$ defined
via (\ref{introduction:gamma0}) at the nodes of the $k$-point Chebyshev extremal grid
on the small interval $[-a_0,a_0]$.
In the second stage, we use $\gamma_0$  to form an initial guess for Newton iterations which converge
to a vector giving the values of the desired  slowly-varying solution $\gamma$ at the  Chebyshev extremal nodes in $[-a_0,a_0]$.
The first-order approximate $\gamma_0$ is used to form an initial guess because, unlike its higher order analogues $\gamma_M$,
it can be readily calculated.
In the final stage, we apply a standard adaptive Chebyshev spectral method, with the computed values  $\gamma$ and
its derivatives at $0$  used as initial conditions,  to compute a piecewise Chebyshev expansion
representing $\gamma$ over the entire solution domain $[a,b]$.  We detail
 each of these three procedures in a subsection below. 


It will be convenient to write the coefficient $q(t)$ of  (\ref{introduction:ode})  in the form $t q_0(t)$, where $q_0$ is smooth and 
positive on $[a,b]$. Our algorithm takes the following as inputs:
\renewcommand\labelitemi{\tiny$\bullet$}
\begin{itemize}
    \item an interval $[a,b]$ containing $0$ and over which the Airy-Kummer phase function $\gamma$ is to be computed;
    \item a sufficiently small subinterval of the form $[-a_0,a_0]$ of $[a,b]$ containing $0$ on which to apply the Newton-Kantorovich method;
    \item a positive integer $k$ which controls the order of the Chebyshev expansions used to represent $\gamma$; 
    \item a positive real number $\epsilon$ specifying the desired precision for the calculations;  
    \item the value of the parameter $\omega$; and
    \item an external subroutine for evaluating $q_0$.
\end{itemize}
It  outputs $(k-1)^{st}$ order piecewise Chebyshev expansions representing the Airy phase function $\gamma$ and its
first two derivatives  on the interval $[a,b]$.
To be entirely clear, a $(k-1)^{st}$ order piecewise Chebyshev expansion on the interval $[a,b]$ is 
an expansion of the form
\begin{equation} 
    \begin{split}
         &\sum_{i=1}^{m-1} \chi_{[\xi_{i-1},\xi_i)}(t) \sum_{j=0}^{k-1} c_{ij} T_j 
\left( \frac{2}{\xi_i-\xi_{i-1}} t - \frac{\xi_i+\xi_{i-1}}{\xi_i-\xi_{i-1}} \right) \\
        &+ \chi_{[\xi_{m-1},\xi_m]}(t) \sum_{j=0}^{k-1} c_{mj} T_j \left( \frac{2}{\xi_m-\xi_{m-1}} t - \frac{\xi_m + \xi_{m-1}}{\xi_m - \xi_{m-1}} \right)
    \end{split}
\label{algorithm:3}
\end{equation}
where $a = \xi_0 < \xi_1 < \cdots < \xi_m = b$ is a partition of $[a,b]$, $\chi_I$ is the characteristic function on the interval $I$
and $T_j$ is the Chebyshev polynomial of degree $j$. We note that the characteristic function of a half-open interval appears in the first line
of (\ref{algorithm:3}), whereas the characteristic function of a closed interval appears in the second line.
 This ensures that exactly one of the characteristic functions appearing in (\ref{algorithm:3}) is nonzero for each point $t$ in $[a,b]$.

Once the piecewise Chebyshev expansions representing $\gamma$ and its first two derivatives have  been constructed, both elements 
of a basis in the space of solutions of (\ref{introduction:ode}) and their first derivatives
can be readily evaluated at any point in $[a,b]$ at a cost which is independent of $\omega$.  Of course, it follows 
that a large class of  initial and boundary value problems
for (\ref{introduction:ode}) can be solved in time independent of $\omega$ and 
the obtained solutions can be evaluated at arbitrary points in
$[a,b]$ in time independent of $\omega$.



%
%
\begin{subsection}{Numerical evaluation of the asymptotic approximation}

To evaluate $\gamma_0$, we first compute the values of the function
\begin{equation} 
 \int_0^{t} \sqrt{|q(s)|} \ ds
=     \int_0^{t} \sqrt{|s|} \sqrt{q_0(s)} \ ds
\label{algorithm:5}
\end{equation}
at the nodes $t_1,\ldots,t_k$ of the Chebyshev extremal grid on $[-a_0,a_0]$.   
Because of the singularity in the integrand, (\ref{algorithm:5})
cannot be evaluated efficiently
using Clenshaw-Curtis or Gauss-Legendre quadrature rules.  We could use a Gauss-Jacobi rule, but
 we prefer another technique based on monomial expansions.  
Although the numerical use of monomial expansions has historically  been viewed with suspicion owing to concerns
regarding numerical stability, it is shown in \cite{ShenSerkh} that  under widely-applicable conditions they are
as stable as an orthogonal polynomial basis when used for interpolation.     

We proceed by numerically calculating the coefficients of a monomial expansion
\begin{equation}
p(t) = \sum_{j=0}^{k-1} c_j \left(\frac{t}{a_0}\right)^{j-1}
\label{algorithm:6}
\end{equation}
representing the smooth function $\sqrt{q_0(t)}$ over $[-a_0,a_0]$.  This is done by solving the  Vandermonde system
which results from enforcing the conditions
\begin{equation}
p(t_i) = \sqrt{q_0(t_i)}, \ \ \ i=1,\ldots,k.
\label{algorithm:7}
\end{equation}
Inserting the expansion (\ref{algorithm:6}) into (\ref{algorithm:5}) and evaluating
the integral yields
\begin{equation}
\begin{aligned}
\sum_{j=0}^{k-1} c_j \int_{0}^{t} \sqrt{|s|}\left(\frac{s}{a_0}\right)^{j-1}\, ds
=
a_0 \sqrt{|t|} \sum_{j=0}^{k-1} \frac{c_j}{j+\frac{1}{2}} \left(\frac{t}{a_0}\right)^{j}.
\end{aligned}
\label{algorithm:8}
\end{equation}
Using (\ref{algorithm:8}), the values of the  function (\ref{algorithm:5}) at the Chebyshev extremal nodes $t_1,\ldots,t_k$ can be easily computed
and, once this has been done,  it is trivial to evaluate the function $\gamma_0$ defined via (\ref{introduction:gamma0})
at those same nodes.  The output of the procedure of this subsection is the vector
\begin{equation}
\boldsymbol{\gamma_0} \approx \left(\begin{array}{cccc}
\gamma_0(t_1) & \gamma_0(t_2) & \cdots & \gamma_0(t_k)
\end{array}\right)^\top
\end{equation}
approximating the values of the first order approximate $\gamma_0$ at the extremal Chebyshev nodes on $[-a_0,a_0]$.

\begin{remark}
It is easy to generalize the procedure of this subsection to a nonsymmetric interval of the form $[a_0,b_0]$ containing $0$.
In this case, for the sake  of numerical stability, it is best to use a monomial expansion of the form
\begin{equation}
p(t) = \sum_{j=0} c_j \left(\frac{2}{b_0-a_0}t -  \frac{b_0+a_0}{b_0-a_0}\right)^{j-1}
\label{algorithm:9}
\end{equation}
to represent $\sqrt{q_0(t)}$ and the integrals which need to be evaluated are of the form
\begin{equation}
\int_{0}^{t} \sqrt{|s|}\left(s -\frac{a_0+b_0}{2}\right)^{j-1}\, ds.
\label{algorithm:10}
\end{equation}
An explicit expression for  (\ref{algorithm:10}) can be written in terms of the Gaussian hypergeometric function
${}_2F_1(a,b;c;z)$, which can be easily evaluated via a three-term linear recurrence relation.  The necessary
formulas can be found, for example, in  Chapter~2 of \cite{HTFI}.
\end{remark}

\end{subsection}

%
%

\begin{subsection}{Computation of $\gamma$ over the interval $[-a_0,a_0]$}

In this subsection, we describe our procedure for calculating the values of 
the desired slowly-varying phase function  $\gamma$ over the small interval $[-a_0,a_0]$.
Our method consists of using the Newton-Kantorovich method to solve
the discretized version $R\left(\boldsymbol{\gamma}\right)=0$ of 
the Airy-Kummer equation.


We use $\boldsymbol{\gamma}$ to denote the current iterate and,
in the first instance, we take it to be equal to the vector $\boldsymbol{\gamma_0}$.
We let $\boldsymbol{q}$ be the vector defined in (\ref{newton:q})
and  then form the vectors
\begin{equation*}
\boldsymbol{\gamma'} = \frac{1}{a_0} \mathscr{D}_k \boldsymbol{\gamma},\ \ \ 
\boldsymbol{\gamma''} = \frac{1}{a_0} \mathscr{D}_k\boldsymbol{\gamma'}\ \ \ \mbox{and}\ \ \ 
\boldsymbol{\gamma'''} = \frac{1}{a_0} \mathscr{D}_k\boldsymbol{\gamma''}
\end{equation*}
that give the values of the first, second and third derivatives of the Chebyshev expansion
represented via $\boldsymbol{\gamma}$ at the Chebyshev nodes.  
Next, we repeatedly perform the following steps:

\begin{enumerate}

\item
Form the vector $\boldsymbol{r} = R\left(\boldsymbol{\gamma}\right)$, where $R$ is
defined via (\ref{newton:R}).

\item
Form the matrix representing the Fr\'echet derivative $D_{\boldsymbol{\gamma}} R$
of the operator $R$ defined in (\ref{newton:R}) at the point $\boldsymbol{\gamma}$.
It is defined via Formula~(\ref{newton:DR}).

\item
Solve the system of linear equations $D_{\boldsymbol{\gamma}} R \, \left(\boldsymbol{h}\right) = - \boldsymbol{r}$.

\item
Compute the quantity $\zeta = \left\|\boldsymbol{h} \right\|_\infty / \left\|\boldsymbol{\gamma} \right\|_\infty$.
%

\item Let $\boldsymbol{\gamma} = \boldsymbol{\gamma}+\boldsymbol{h}$,
$\boldsymbol{\gamma'} = \frac{1}{a_0} \mathscr{D}_k \boldsymbol{\gamma}$,
$\boldsymbol{\gamma''}  = \frac{1}{a_0} \mathscr{D}_k\boldsymbol{\gamma'}$
and
$\boldsymbol{\gamma'''}  = \frac{1}{a_0} \mathscr{D}_k\boldsymbol{\gamma''}$.

\item If $\zeta > \epsilon$, goto step 1.  Otherwise,  the procedure  terminates.

\end{enumerate}

Upon termination of the above procedure, we have the values of Chebyshev expansions
representing the desired slowly-varying Airy phase function $\gamma$ 
and its first three derivatives on the interval $[-a_0,a_0]$.  We now 
use these expansions to calculate approximations of $\gamma(0)$, $\gamma'(0)$  and $\gamma''(0)$,  
and these quantities comprise the output of this stage of our algorithm.

\end{subsection}

%
%

\begin{subsection}{Extension of $\gamma$ to $[a,b]$}

We now use a  standard  adaptive Chebyshev spectral method to 
solve the  Airy-Kummer equation over the entire interval $[a,b]$.  
We impose the conditions that the obtained solution and its first two derivatives agree
with the values of  $\gamma(0)$, $\gamma'(0)$ and $\gamma''(0)$ computed in the preceding 
stage of our algorithm.

We describe the solver's operation in the case of the more general problem
\begin{equation}
    \begin{cases} 
        \boldsymbol{y}'(t) &= F(t,\boldsymbol{y}(t)), \qquad a < t < b, \\
       \boldsymbol{y}(0) &= \boldsymbol{v}
    \end{cases}
\label{extension:1}
\end{equation}
where $F: \mathbb{R}\times \mathbb{R}^n \to \mathbb{R}^n$ is smooth and $\boldsymbol{v} \in \mathbb{R}^n$.  
Obviously, the initial value problem for the the  Airy-Kummer equation we seek to solve can be put into the form (\ref{extension:1}).
The spectral solver outputs $n$ piecewise $(k-1)^{st}$ order Chebyshev expansions, 
one for each of the components $y_i(t)$ of the solution $\boldsymbol{y}$ of (\ref{extension:1}). 

The solver proceeds in two stages. In the first, it constructs the solution over the interval $[0,b]$. 
During this stage, two lists of subintervals of $[0,b]$ are maintained:
one consisting of what we term  ``accepted subintervals" and the other of subintervals which have yet to be 
processed. A subinterval is accepted if the solution is deemed to be adequately represented by a $(k-1)^{st}$ order 
Chebyshev expansion on that subinterval. Initially, the list of accepted subintervals is empty and the list of subintervals to 
process contains the single interval $[0,b]$.  The solver proceeds as follows until the list of subintervals to process is empty:
\begin{enumerate}
    \item Find, in the list of subintervals to process, the interval $[c,d]$ such that $c$ is as small as possible and remove this subinterval from the list
    \item Solve the initial value problem
    \begin{equation} 
      \left\{
        \begin{split}
            \boldsymbol{u}'(t) &= F(t,\boldsymbol{u}(t)), \qquad c < t < d \\
            \boldsymbol{u}(c) &= \boldsymbol{w}
        \end{split}\right.
        \label{extension:2}
    \end{equation}
    If $[c,d] = [0,b]$, then we take $\boldsymbol{w} = \boldsymbol{v}$. Otherwise, the value of the solutions at the point $c$ 
has already been approximated, and we use that estimate for $\boldsymbol{w}$ in (\ref{extension:2}). If the problem is linear, a 
straightforward Chebyshev integral equation method is used to solve (\ref{extension:2}). Otherwise, the trapezoidal method is first used to 
produce an initial approximation $\boldsymbol{y_0}$ of the solution and then Newton's method is applied to refine it. The linearized problems 
are solved using a Chebyshev integral equation method. In any event, the result is a set of $(k-1)^{st}$ order Chebyshev expansions
    \begin{equation} 
        u_i(t) \approx \sum_{j=0}^{k-1} c_{ij} T_j \left( \frac{2}{d-c} t - \frac{d+c}{d-c} \right), \qquad i = 1,\dots,n
        \label{extension:3}
    \end{equation}
    which purportedly approximate the components $u_1, \dots, u_n$ of the solution of (\ref{extension:2}). 
    \item Compute the quantities 
    \begin{equation} \label{extension:4}
        \frac{\sqrt{\sum_{j=\lfloor \frac{k}{2} \rfloor+1}^{k-1} |c_{ij}|^2 }}{\sqrt{\sum_{j=0}^{k-1} |c_{ij}|^2}}, \qquad i = 1,\dots,n
    \end{equation}
    where $c_{ij}$ are the coefficients in the expansions (\ref{extension:3}). If any of the resulting values is larger than $\epsilon$, then 
split the subinterval into two halves $\left( c,\frac{c+d}{2} \right)$ and $\left(\frac{c+d}{2},d\right)$ and place them on the list of 
subintervals to process. Otherwise,  place the subinterval $(c,d)$ on the list of accepted subintervals.
\end{enumerate}
At the conclusion of this stage, we have $(k-1)^{st}$ order piecewise Chebyshev expansions representing
 each  component of the solution over the interval $[0,b]$, with the list of accepted subintervals determining the partition of $[0,b]$
associated with the piecewise expansions.

In its second stage, an analogous procedure is used to construct piecewise Chebyshev expansions
representing the  solution over the interval $[a,0]$.  In each step, instead of choosing
the unprocessed interval $[c,d]$ such that $c$ is as small as possible and solving an initial value problem
over $[c,d]$, a terminal value problem of the form
    \begin{equation} 
\left\{
        \begin{split}
            \boldsymbol{u}'(t) &= F(t,\boldsymbol{u}(t)), \qquad c < t < d \\
            \boldsymbol{u}(d) &= \boldsymbol{w},
        \end{split}\right.
        \label{extension:5}
    \end{equation}
where $[c,d]$ is the  unprocessed interval such that $d$ is a large as possible, is solved.
At the conclusion of this second stage, we have a $(k-1)^{st}$ order piecewise Chebyshev expansion
representing each component of the solution over the interval $[a,0]$, with the list of accepted subintervals determining the partition of $[a,0]$
associated with the piecewise expansions.  Obviously, amalgamating the piecewise expansions
of the $u_j$ produced during the two stages gives us the desired piecewise Chebyshev expansions of the components
of the solution of (\ref{extension:1}) over the entire interval $[a,b]$.

\end{subsection}

\end{section}

\begin{section}{Numerical Experiments}
\label{section:experiments}

In this section, we present the results of numerical experiments conducted to illustrate 
the properties of our method.   We implemented our algorithm in Fortran and compiled our 
code with  version  14.2.1 of the  GNU Fortran compiler.    All experiments were performed on a desktop 
computer equipped with an AMD 9950X processor and 64GB of RAM.  This processor  has 16 cores, but
only one was utilized in our experiments.  In all of our experiments, we took the parameter
$k$ controlling the order of the piecewise Chebyshev expansions used to 
 represent Airy phase functions to be $16$, and we set
precision parameter $\epsilon$ to be $10^{-13}$.  To account for the vagaries of modern 
computing environments, all reported times were obtained  
by averaging the cost of each calculation over 100 runs.

For the most part, we measured the accuracy of our method by using it to solve
various initial and boundary value problems for second order equations of the form (\ref{introduction:ode}).
Because the condition numbers of such problems grow with the parameter $\omega$,
the accuracy of any  solver will deteriorate with increasing $\omega$.
In the case of our algorithm, the mechanism by which accuracy is lost is well understood.
It computes  the Airy phase functions themselves to high precision,
but the magnitude of the phase functions  increases with the parameter $\omega$,
with the consequence that accuracy is lost when the Airy functions are evaluated at large arguments
in order to calculate the solutions of the original differential equation.
Our algorithm does, however, compute solutions of the original differential
equation with accuracy on the order of that predicted by the condition number of the problem being solved.

There is one experiment in which we measured the accuracy 
of the Airy phase functions produced by our algorithm directly.
In the experiment of Subsection~\ref{section:experiments:3},
we constructed Airy phase functions representing associated Legendre functions
of various degrees and orders  using the algorithm of this paper and then
compared them with Airy phase functions produced via another technique in order
to show our algorithm computes Airy phase functions with  high relative accuracy.

%
%

\begin{subsection}{Initial value problems}
\label{section:experiments:1}

\begin{figure}[t!!!!]
\hfil
\includegraphics[width=\textwidth]{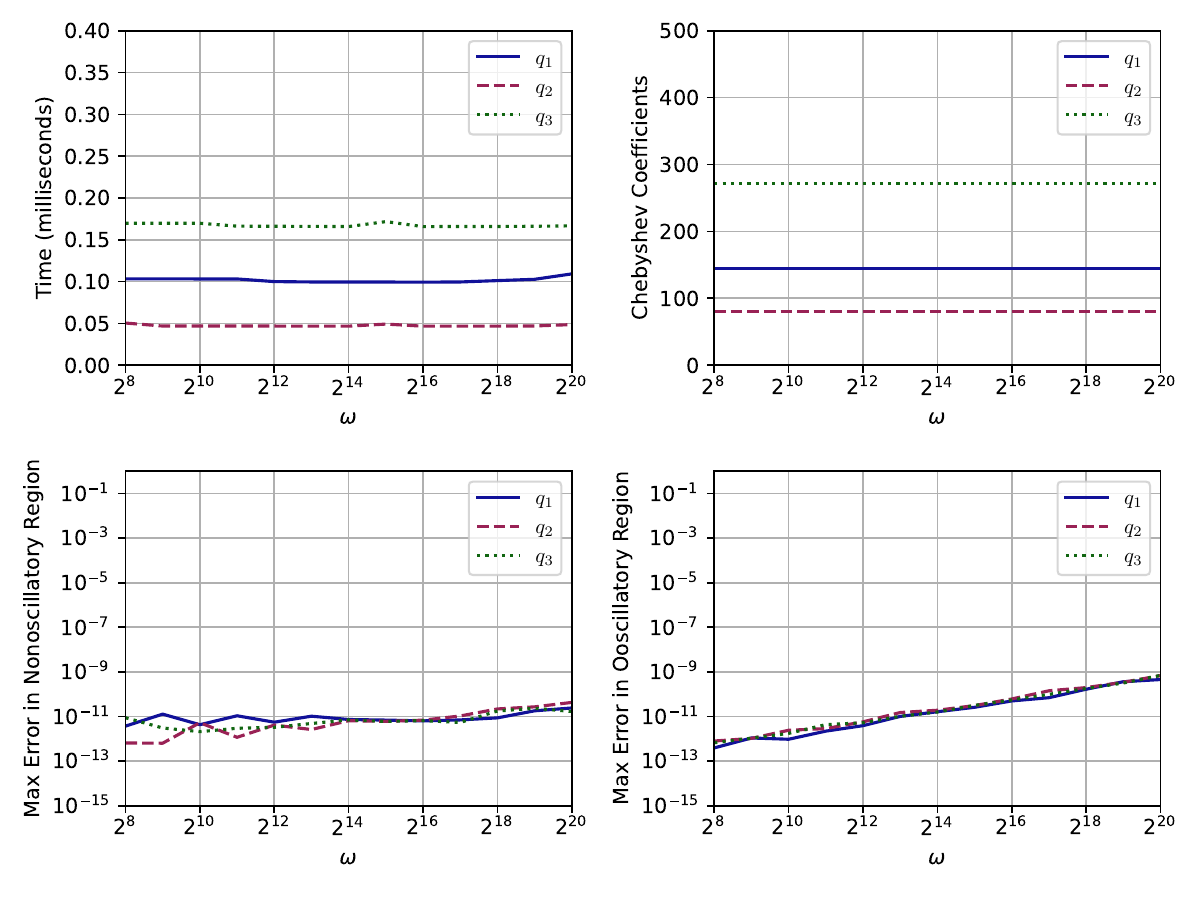}
\hfil
\vskip -2em
\caption{
The results of the experiment of Section~\ref{section:experiments:1}
in which a collection of initial value problems were solved using our method.
The upper-left plot gives the time, in milliseconds, required to compute
each Airy phase function.  The plot in the upper right gives the number
of coefficients in the piecewise Chebyshev expansions used to represent
the Airy phase functions. The plot in the lower left gives the
relative accuracy of the obtained solutions of the initial value problems
in the nonoscillatory region, and the plot in the lower right
reports the absolute accuracy of the obtained solution in the oscillatory region.
}
\label{experiments1:figure1}
\end{figure}

In our first experiment, we used the algorithm of this paper to solve
initial value problems of the form
\begin{equation}
\left\{
\begin{aligned}
&y''(t) + \omega^2 q(t) y(t) = 0, \ \ \ -5 <t < 5,\\
&y(0) = 1,\ \ y'(0) = 0
\end{aligned}
\right.
\label{experiments:1:ivp}
\end{equation}
for various values of $\omega$ and choices of the coefficient $q(t)$.  More explicitly, 
for each of the coefficients
\begin{equation}
\begin{aligned}
q_1(t) = t + t^3, \ \ \ q_2(t) = -1 + (1-t)\exp(t)\ \ \ \mbox{and}
\ \ \ q_3(t) = t + \frac{\sin(3t)}{3}
\end{aligned}
\end{equation}
and each $\omega=2^8,\, 2^9,\, 2^{10},\, \ldots,\, 2^{20}$, we used the algorithm of 
Section~\ref{section:algorithm} to construct a slowly-varying Airy phase function $\gamma$ 
representing the solutions of (\ref{experiments:1:ivp}).  We then used each of these Airy phase functions to 
calculate the solution of (\ref{experiments:1:ivp}) at a collection of
evaluation points in the interval $(-5,5)$ and compared the obtained values with a reference solution
constructed via an adaptive Chebyshev spectral method.

\begin{figure}[t!!!!!!]
\hfil
\includegraphics[width=\textwidth]{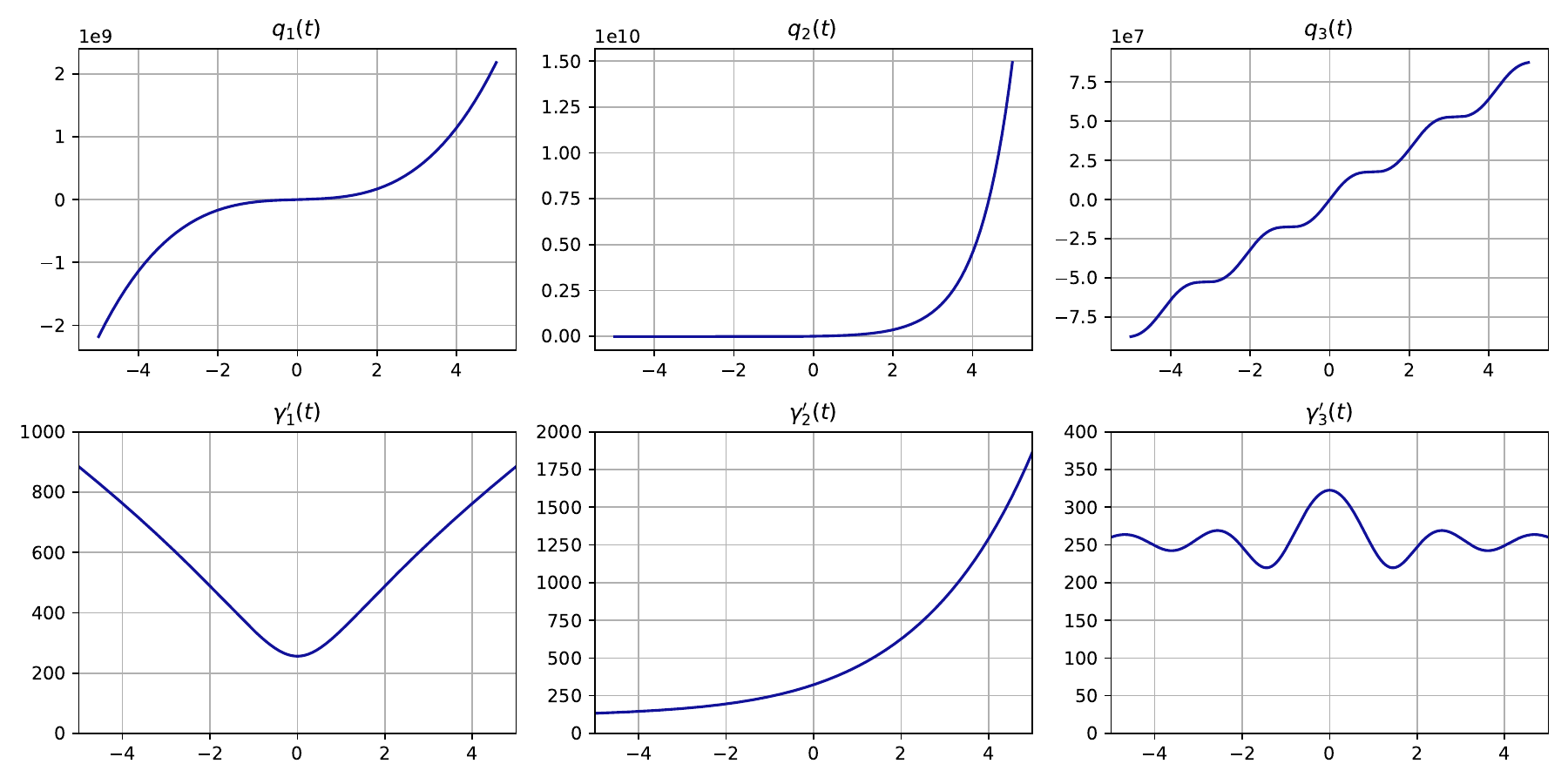}
\hfil
\vskip -1em
\caption{
Plots of the coefficients $q_1(t)$, $q_2(t)$ and $q_3(t)$ considered in the experiment of Section~\ref{section:experiments:1}
when the parameter $\omega$ is equal to $2^{12}$,
as well as plots of the derivatives of the corresponding Airy phase functions $\gamma_1'(t)$, $\gamma_2'(t)$ 
and $\gamma_3'(t)$.
}
\label{experiments1:figure2}
\end{figure}


Since each of the coefficients we consider is positive in the interval $(0,5)$ and negative 
in the interval $(-5,0)$, the solutions of (\ref{experiments:1:ivp}) oscillate  in $(0,5)$
and are nonoscillatory in $(-5,0)$.  Moreover, the solutions increase very rapidly as the argument $t$
decreases from $0$ to $-5$.  Indeed, in each case we considered, the solutions of (\ref{experiments:1:ivp}) 
were too large to  represent via double precision numbers on much of the interval $(-5,0)$.   Accordingly, 
for each problem, we   calculated the largest  relative error in the obtained solution 
at  1,000  equispaced evaluation points in the subinterval $\left[\widetilde{a},0\right]$, 
where  $\widetilde{a}$ was chosen so that $\gamma\left(\widetilde{a}\right) = -100$,
and we used this quantity as our measure of the  accuracy of the solution of (\ref{experiments:1:ivp})
obtained via our method in the nonoscillatory regime.   Since
\begin{equation}
\Ai(-100) \approx 4.669498035610554\times 10^{-291}
\ \ \ \mbox{and}\ \ \ 
\Bi(-100) \approx 1.070779073708091 \times 10^{289},
\end{equation}
$\widetilde{a}$ is close to the point at which the solutions become too large to represent using double
precision numbers.    We note that the Airy phase functions themselves are computed 
accurately over the entire interval $[-5,5]$ via our method.
We also used the Airy phase functions to evaluate each solution 
at  $1,000$ equispaced evaluation points in the interval $[0,5]$
and compared those values with the reference solution.  Since the solutions are oscillatory
here, we measured absolute rather than relative errors at these evaluation points.

Figure~\ref{experiments1:figure1} gives the results of this experiment.
We observe that the running time of our algorithm and the number of Chebyshev coefficients
needed to represent the Airy phase functions are both independent of the parameter $\omega$.
The errors in the obtained solutions of the initial value problems, on the other hand, increase as
$\omega$ grows.   This is to be expected since the condition number of the 
each of the problems we considered increases with $\omega$ and any numerical algorithm for 
solving these problems will lose accuracy as $\omega$ increases.
Plots of the coefficients  $q_1(t)$, $q_2(t)$ and $q_3(t)$  and the 
derivatives of the corresponding Airy phase functions $\gamma_1'(t)$,
$\gamma_2'(t)$ and  $\gamma_3'(t)$ when $\omega=2^{12}$ can  be found in Figure~\ref{experiments1:figure2}.

\end{subsection}

%
%

\begin{subsection}{Boundary value problems}
\label{section:experiments:2}

In our next experiment, we used the algorithm of this paper to solve
the boundary value problem
\begin{equation}
\left\{
\begin{aligned}
&y''(t) + \omega^2 q(t, \omega) y(t) = 0, \ \ \ 0<t<3,\\
&y(0) = 1, \ \ y(3) = 1
\end{aligned}
\right.
\label{experiments:2:bvp}
\end{equation}
for each $\omega=2^8,\, 2^9,\, 2^{10},\, \ldots,\, 2^{20}$ and the following choices of $q$:
\begin{equation}
\begin{aligned}
q_1(t,\omega) = t+t^3, \ \ \ 
q_2(t,\omega) = \sin(t) + 2 \sin\left(\frac{t}{4}\right)^2\ \ \ \mbox{and} \ \ \ 
q_3(t,\omega) = \frac{t \left(\cos ^2(3 t) \sin ^2(\omega )+2\right)}{t^2 \cos ^2(\omega)+1}.
\end{aligned}
\end{equation}
We tested the accuracy of each solution by comparing its value at 1,000 equispaced evaluation points in the 
solution domain $(0,3)$ to a reference solution constructed via an adaptive Chebyshev spectral method.    
Since the solutions are oscillatory in the interval $(0,3)$, we measured absolute rather than relative errors.

The results are given in  Figure~\ref{experiments2:figure1}.    As with the experiments of the 
preceding section, the errors in the obtained  solutions increase with $\omega$.    But again, this is to be
expected since the condition numbers of these boundary value problems grow with $\omega$ and a similar
loss of accuracy will be experienced by any numerical method.  We also note that the errors exhibit greater 
variability than in the previous experiment.  This is due to greater variance 
in the condition numbers of the boundary value problems we considered.
Finally, we observe that the time required to compute the Airy phase function
in the case of the coefficient $q_3$ varies noticeably with $\omega$ because, unlike the
other coefficients, $q_3$ depends on $\omega$.

\begin{figure}[h!!!!!!!!!!!!!!!!!!!!]
\hfil
\includegraphics[width=\textwidth]{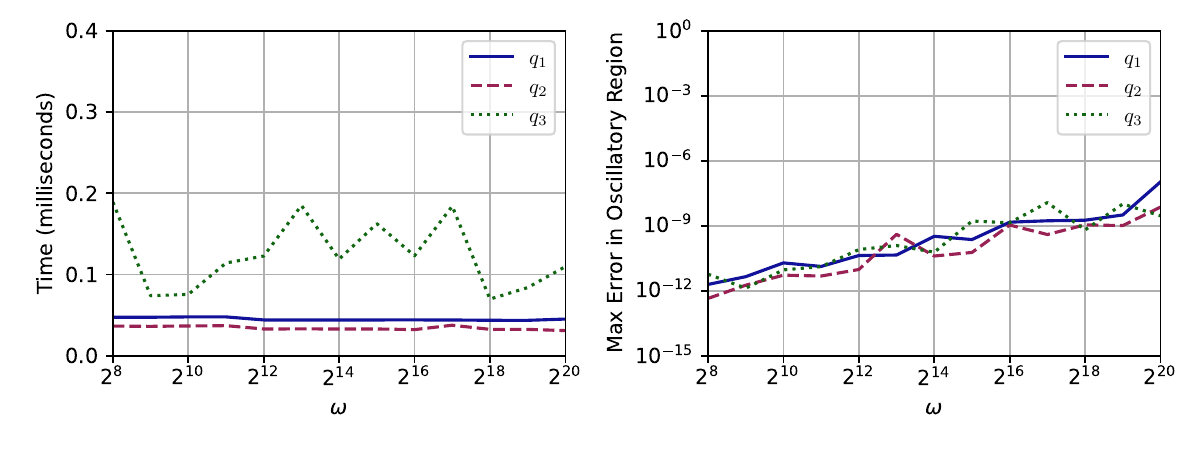}
\hfil
\vskip -2em
\caption{
The results of the experiment of Section~\ref{section:experiments:2}
in which a collection of boundary value problems were solved using our method.
The plot on the left gives the time, in milliseconds, required to compute
each Airy phase function, while the plot on the right gives the
accuracy of the obtained solutions of the boundary  value problems.
}
\label{experiments2:figure1}
\end{figure}
\end{subsection}

%
%
%
%

\begin{subsection}{The associated Legendre differential equation}
\label{section:experiments:3}

\begin{figure}[t!]
\hfil
\includegraphics[width=\textwidth]{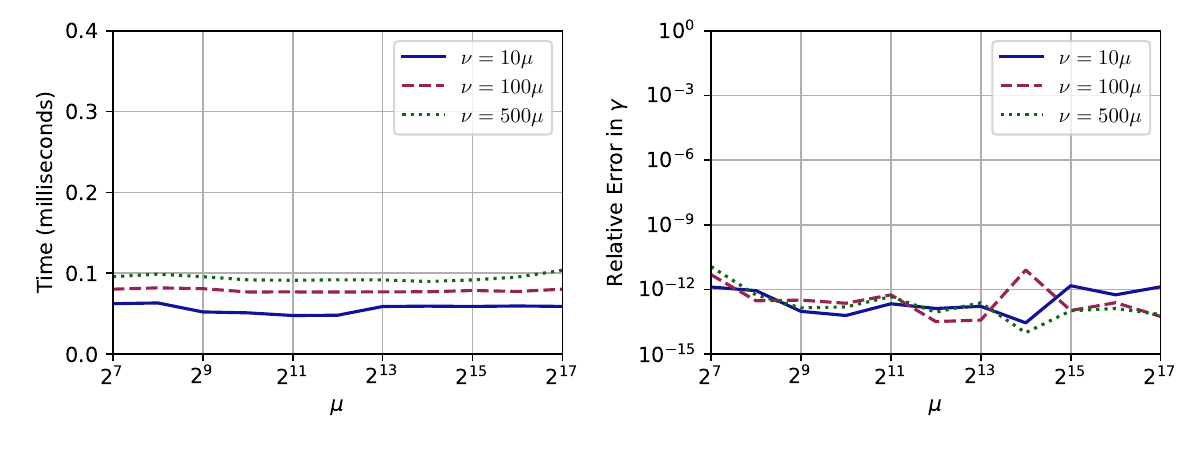}
\hfil
\vskip -2em
\caption{
Some of the results of the experiment of Section~\ref{section:experiments:3},
which concern the associated Legendre functions.
The plot on the left gives the time required to compute each Airy phase function,
while the plot on the right gives the maximum relative error in the 
 Airy phase functions calculated using the method of this paper.
}
\label{experiments3:figure1}
\end{figure}

The standard solutions of the associated Legendre differential equation
\begin{equation}
(1-t^2) y''(t) -2t y'(t) + \left(\nu(\nu+1) -\frac{\mu^2}{1-t^2}\right) y(t) = 0
\label{experiments:3:lode}
\end{equation}
on the interval $(-1,1)$ are the Ferrers functions of the   first and second kinds $P_\nu^\mu$  and $Q_\nu^\mu$.
Definitions of them can be found in Section~14.3 of \cite{DLMF} or Section~5.15 of \cite{Olver}.
   Equation~(\ref{experiments:3:lode}) has singular points at $\pm 1$, and, as long as $\nu > \mu$,
it has  two turning points in the interval $(-1,1)$.  In order to put this equation into a form suitable for our algorithm, 
we introduce the  change of variables $t = \tanh\left(x + \xi_\nu^\mu \right)$, where  
\begin{equation}
\xi_\nu^\mu = \arccosh\left(\frac{\sqrt{\nu(\nu+1)}}{\mu}\right),
\label{experiments:3:x0}
\end{equation}
which yields the new equation
\begin{equation}
y''(x) + \left(-\mu^2 + \nu(\nu+1) \sech^2(x-\xi_\nu^\mu)\right) y(x) = 0.
\label{experiments:3:transformed}
\end{equation}
The coefficient in (\ref{experiments:3:transformed}) is  smooth and while it has two turning points 
located at $0$ and $2 \xi_\nu^\mu$,  because of the symmetries of the associated Legendre functions, it suffices 
to consider it on the interval $\left(-\infty,\xi_\nu^\mu\right)$ which contains only the turning point at $0$.


For each $\mu=10^7,\, 10^8,\, \ldots,\, 10^{17}$ and each $\nu=10\mu,\, 100\mu,\, 500\mu$, we used the algorithm
of this paper to construct  an Airy phase function $\gamma_\nu^\mu$ 
representing the solutions of (\ref{experiments:3:transformed})
over the interval $\left[a_\nu^\mu,\xi^\mu_\nu\right]$, where $\xi_\nu^\mu$ is as in (\ref{experiments:3:x0}) and $a_\nu^\mu$ 
is chosen such that $\gamma_\nu^\mu\left(a_\nu^\mu\right) = -15$.
We note that 
\begin{equation}
\Bi(-15) \approx  3.364489547667594\times 10^{16}
\ \ \ \mbox{and}\ \ \ 
\Ai(-15) \approx 3.837296156948168\times10^{-18},
\end{equation}
so that the solutions of (\ref{experiments:3:transformed}) are either of very large or very
small magnitude at $a_\nu^\mu$.  

We first tested the accuracy of each $\gamma_\nu^\mu$ by 
using it to evaluate a solution of (\ref{experiments:3:transformed}) at $1,000$ equispaced 
points in the solution domain and recording the largest  observed relative error.  
Because the  Ferrers functions themselves  are normalized such that
 $P_\nu^\mu$ and $Q_\nu^\mu$ become  astronomically large  on portions of $(-1,1)$, even for relatively 
small values of  $\nu$ and $\mu$, we choose to evaluate a different solution instead.
More explicitly, we considered the solution
\begin{equation}
F_\nu^\mu(x) = \widetilde{Q}_\nu^\mu(x) + i \widetilde{P}_\nu^\mu(x),
\label{experiments:3:F}
\end{equation}
where $\widetilde{Q}_\nu^\mu$  and $\widetilde{P}_\nu^\mu$  are given via the formulas
\begin{equation}
\begin{aligned}
\widetilde{Q}_\nu^\mu(x) =  
\sqrt{\frac{2}{\pi}\, \frac{\Gamma\left(1+\nu+\mu\right)}
{\Gamma\left(1+\nu-\mu\right)}} Q_\nu^{-\mu}(\tanh(x))
\ \ \ \mbox{and}\ \ \ 
\widetilde{P}_\nu^\mu(x) 
= \sqrt{\frac{\pi}{2}\, \frac{\Gamma\left(1+\nu+\mu\right)}
{\Gamma\left(1+n-m\right)}} P_\nu^{\mu}(\tanh(x)).
\end{aligned}
\label{experiments:3:normalizedpq}
\end{equation}
The Wronskian of this pair is $1$, which ensures that $F_\nu^\mu$ is normalized in a reasonable way.
Moreover, $\widetilde{P}_\nu^\mu$ and $\widetilde{Q}_\nu^\mu$ determine a slowly-varying trigonometric
phase function $\alpha_\nu^\mu$ for (\ref{experiments:3:transformed}) through the relations
\begin{equation}
\widetilde{Q}_\nu^\mu(x) = \frac{\cos\left(\alpha_\nu^\mu(x)\right)}{\sqrt{\frac{d}{dx} \alpha_\nu^\mu(x)}}, \ \ \
\widetilde{P}_\nu^\mu(x) = \frac{\sin\left(\alpha_\nu^\mu(x)\right)}{\sqrt{\frac{d}{dx} \alpha_\nu^\mu(x)}}
\ \ \ \mbox{and}\ \ \
\lim_{x\to -\infty} \alpha_\nu^\mu (x) = 0
\label{experiments:3:pqalpha}
\end{equation}
It is beyond the scope of this paper to show that $\alpha_\nu^\mu$ is, in fact, slowly varying, but we refer
the interested reader to \cite{BremerPhase2}.  Since the logarithmic derivative of $F_\nu^\mu$ is the derivative of $\alpha_\nu^\mu$,
its  condition number of evaluation is slowly varying.

\begin{figure}[t!]
\hfil
\includegraphics[width=\textwidth]{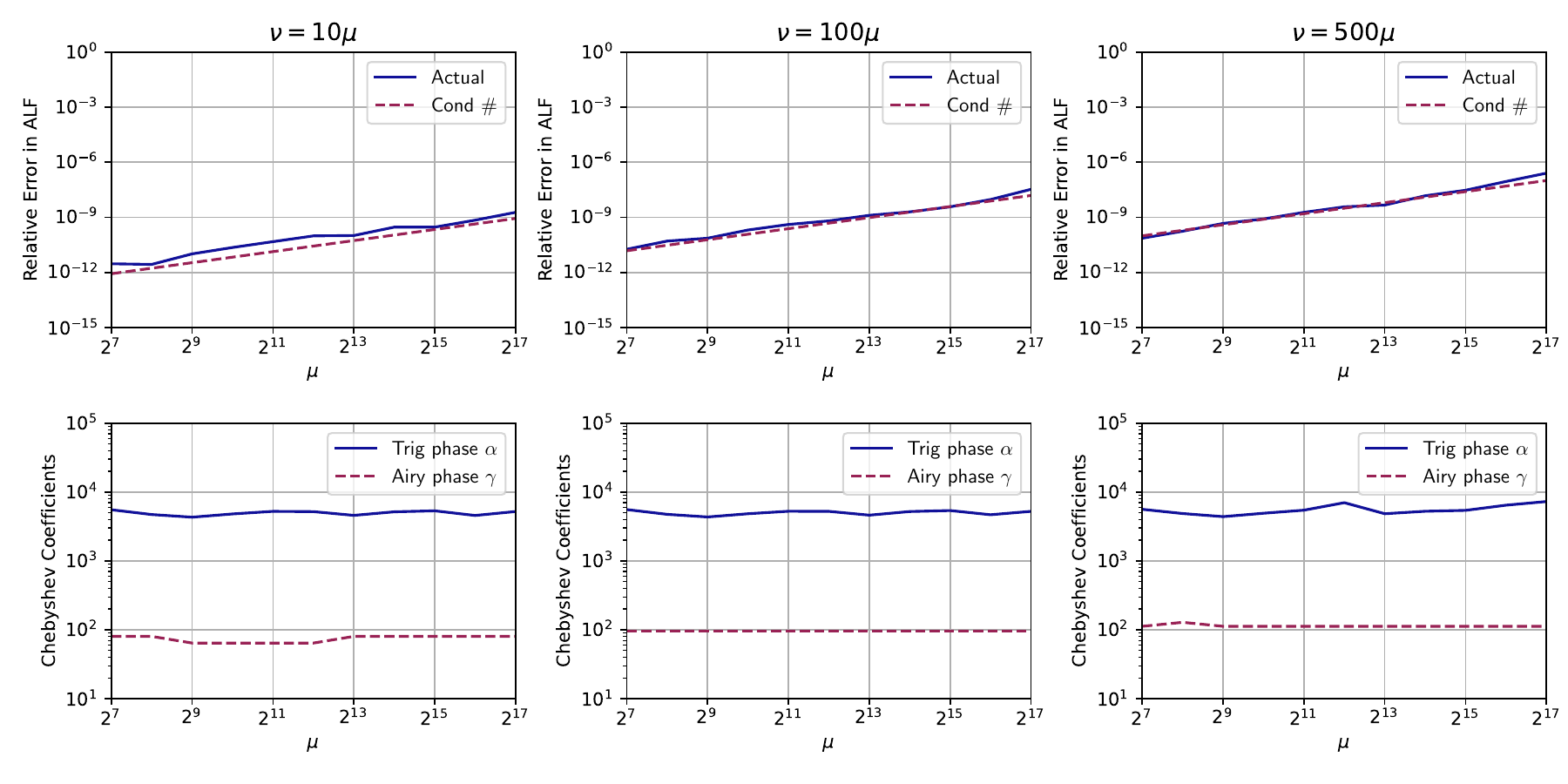}
\hfil
\vskip -2em
\caption{
Some of the  results of the experiment of Section~\ref{section:experiments:3},
which concerns the associated Legendre functions.  The plots on the top
row compare the relative error in the evaluation
of the function $F_\nu^\mu$ defined via (\ref{experiments:3:F}) 
with the accuracy predicted by its condition number of evaluation.
The plots on the bottom row report the number of coefficients
in the piecewise Chebyshev expansions of the Airy phase functions $\gamma_\nu^\mu$
and the trigonometric  phase functions $\alpha_\nu^\mu$.
}
\label{experiments3:figure2}
\end{figure}

Next, for each pair of values of $\nu$ and $\mu$ considered, 
we used an alternate approach to construct a second Airy phase function  $\widetilde{\gamma}_\nu^\mu$
and measured the relative accuracy of $\gamma_\nu^\mu$ by comparing the two phase functions
at 1,000 equispaced points on the interval $\left[a_\nu^\mu,\xi^\mu_\nu\right]$.
To construct $\widetilde{\gamma}_\nu^\mu$,
we first used  the algorithm of \cite{StojimirovicBremer} to calculate
the slowly-varying trigonometric phase function  $\alpha_\nu^\mu$ defined
via (\ref{experiments:3:pqalpha}) and then composed it with the inverse
of a slowly-varying trigonometric phase function $\aiphase$ for Airy's equation;
that is, we let
$\widetilde{\gamma}_\nu^\mu(t) = \aiphase^{-1}\left(\alpha_\nu^\mu(t)\right)$.
The phase function $\aiphase$ is determined via the requirements
\begin{equation}
\Bi(x) = \frac{\cos\left(\aiphase(x)\right)}{\sqrt{\aiphase'(x)}},\ \ \ 
\Ai(x) = \frac{\sin\left(\aiphase(x)\right)}{\sqrt{\aiphase'(x)}}
\ \ \ \mbox{and}\ \ \
\lim_{x\to -\infty} \aiphase (x) = 0.
\end{equation}
That the composition $\widetilde{\gamma}_\nu^\mu$ is an Airy phase function for (\ref{experiments:3:transformed}) 
follows from the formulas
\begin{equation}
\begin{aligned}
\frac{\Bi\left(\widetilde{\gamma}_\nu^\mu(x) \right)}{\sqrt{\frac{d}{dx} \widetilde{\gamma}_\nu^\mu(x)}}
&=
\frac{\cos\left(\aiphase\left(\widetilde{\gamma}_\nu^\mu(x)\right) \right)}{\sqrt{\alpha_{\mbox{\tiny ai}}'\left(\gamma_\nu^\mu(x)\right)
\frac{d}{dx} \widetilde{\gamma}_\nu^\mu(x)}} 
= 
\frac{\cos\left(\alpha_\nu^\mu(x)\right)}{\sqrt{\frac{d}{dx} \alpha_\nu^\mu(x)}}
= \widetilde{Q}_\nu^\mu(x) \ \ \ \mbox{and}
\\[1em]
\frac{\Ai\left(\widetilde{\gamma}_\nu^\mu(x) \right)}{\sqrt{\frac{d}{dx} \widetilde{\gamma}_\nu^\mu(x)}}
&=
\frac{\sin\left(\aiphase\left(\widetilde{\gamma}_\nu^\mu(x)\right) \right)}{\sqrt{\aiphase'\left(\gamma_\nu^\mu(x)\right)
\frac{d}{dx} \widetilde{\gamma}_\nu^\mu(x)}}=
\frac{\sin\left(\alpha_\nu^\mu(x)\right)}{\sqrt{\frac{d}{dx} \alpha_\nu^\mu(x)}}
=
\widetilde{P}_\nu^\mu(x).
\end{aligned}
\end{equation}

\begin{figure}[b!!!!!!!!]
\hfil
\includegraphics[width=\textwidth]{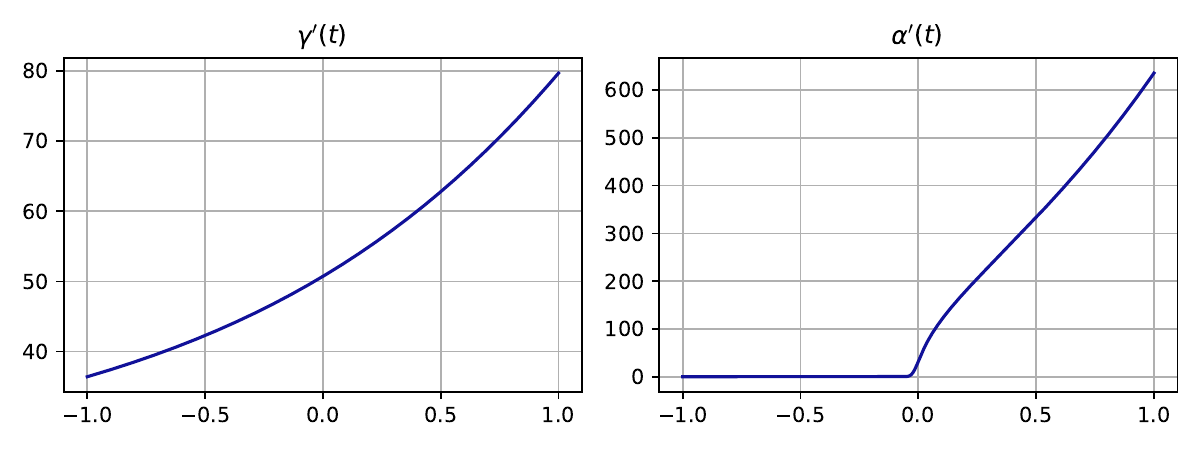}
\hfil
\vskip -2em
\caption{
Graphs of the derivatives of the  Airy phase function $\gamma_\nu^\mu$ (left) and the trigonometric phase  $\alpha_\nu^\mu$ (right)
over the interval $[-1,1]$ in the case $\mu=256$ and $\nu=2560$.       The derivative of the 
trigonometric phase function is significantly more expensive to represent via polynomial expansions than the Airy phase function, 
which is extremely benign.
}
\label{experiments3:figure3}
\end{figure}

The results of this experiment are given in Figures~\ref{experiments3:figure1} and \ref{experiments3:figure2}.
We note that, in addition to the maximum relative error observed while evaluating
$F^\mu_\nu$ using our algorithm,  the plots on the top row of Figure~\ref{experiments3:figure2} report the 
maximum relative accuracy that is expected given the condition number of evaluation of $F_\nu^\mu$,
which is, of course,
\begin{equation}
\kappa_\nu^\mu =  \max_{x_1,\ldots,x_{1000}} \left|x_j \frac{\frac{d}{dx} F_\nu^\mu(x_j)}{ F_\nu^\mu(x_j)}\right|\, \epsilon_0,
\end{equation}
where $\epsilon_0\approx 2.220446049250313\times 10^{-16}$ is machine zero for the IEEE double precision number system.
The plots on the  bottom row of Figure~\ref{experiments3:figure2} compare the number of coefficients
in the piecewise expansions of the trigonometric phase functions $\alpha_\nu^\mu$
and of the Airy phase functions $\gamma_\nu^\mu$.
While both functions are represented at a cost which is independent of $\mu$, the trigonometric
phase function is several orders of magnitude more expensive to represent than  the Airy phase function.
Figure~\ref{experiments3:figure3}, which contains plots of the derivatives of the 
trigonometric phase function $\alpha_\nu^\mu$ and of the Airy phase function $\gamma_\nu^\mu$ 
in the case $\nu = 2560$ and $\mu = 256$, makes clear why this is the case.
In particular, it shows that the derivative of $\alpha_\nu^\mu$ exhibits complicated behavior
near the turning point, while $\gamma_\nu^\mu$ is extremely benign throughout the interval.

\end{subsection}

\end{section}

\begin{section}{Conclusion}
\label{section:conclusion}

We have given a proof of the existence of slowly-varying Airy phase functions and described a numerical
method for rapidly computing them.  Using our algorithm, a large class of second
order linear ordinary differential equations of the form (\ref{introduction:ode})
can be solved to high accuracy in time independent of the parameter $\omega$.  This class
includes many differential equations defining  widely-used special functions,
as well as many equations with applications in physics and chemistry.

With some modification, our numerical method extends to the case of equations of the form
\begin{equation}
y''(t) + \omega^2 t^\sigma q(t) y(t) = 0,
\label{conclusion:1}
\end{equation}
where $q(t)\sim 1$ as $t\to 0$, $\sigma >-2$ and $\omega$ is large.    However, the analysis of this paper
fails when $\sigma$ is not one of the special values  $-1$, $0$ or $1$.  
The authors will discuss the necessary modifications to our numerical algorithm and an alternative 
approach to the procedure of Section~\ref{section:analysis:formal} in a future work.

In the experiment of Subsection~\ref{section:experiments:3}, a second procedure for constructing 
slowly-varying Airy phase functions was introduced.  Namely, the inverse of a 
trigonometric phase function for Airy's equation was composed with a trigonometric
phase function for the associated Legendre equation in order to form a slowly-varying
Airy phase function for the associated Legendre equation.    Methods of this 
type might be useful for various numerical computations as well
as for the derivation of asymptotic approximations for the solutions of ordinary differential equations,
and this approach warrants further investigation.


Finally, we note that the efficient representation of solutions of    second
order linear ordinary differential equations via generalized phase functions
has many applications to the rapid evaluation of special functions and their zeros,
and to the rapid application of the related Sturm-Liouville transforms.    The authors
also plan to explore these topics in future works.

\end{section}

\begin{section}{Acknowledgments}
JB  was supported in part by NSERC Discovery grant  RGPIN-2021-02613.
\end{section}


\bibliographystyle{plain}
\bibliography{airyphase.bib}

\end{document}